\documentclass[]{interact}

\usepackage{epstopdf}

\usepackage[longnamesfirst,sort]{natbib}
\bibpunct[, ]{(}{)}{;}{a}{,}{,}

\theoremstyle{plain}

\theoremstyle{definition}

\theoremstyle{remark}

\usepackage{natbib}
\usepackage{amssymb,amsmath,fancyhdr,supertabular,longtable,hhline}
\usepackage{caption}
\usepackage{enumitem}
\usepackage{subcaption}
\usepackage[colorlinks, hyperindex, plainpages=false, linkcolor=blue, urlcolor=blue, pdfpagelabels]{hyperref}
\usepackage[all]{hypcap}
\usepackage{pgf,tikz}
\usetikzlibrary{arrows}
\usepackage{txfonts,pxfonts}
\usepackage{tikz}
\usetikzlibrary{fit,shapes,patterns,arrows,calc}
\usepackage{tkz-euclide}
\usepackage{dirtytalk}
\usetikzlibrary{matrix}
\usepackage[most]{tcolorbox}
\definecolor{mycolor}{rgb}{0.85, 0.85, 0.85}
\usepackage{pythonhighlight}
\usepackage{enumitem}
\usetikzlibrary{lindenmayersystems}
\usetikzlibrary[shadings]
\usetikzlibrary{decorations.fractals}
\usepackage{logicpuzzle}
\usepackage{enumitem}
\newcommand{\mat}[4]{\left[\begin{smallmatrix}
		#1 & #2 \\
		#3 & #4 \\
	\end{smallmatrix}\right]}

\newcommand{\mbs}{\begin{adjustbox}{minipage=\linewidth,frame}}
	\newcommand {\mbe}{\end{adjustbox}}

\newenvironment{shidoku}[1][]{%
	\begin{logicpuzzle}[rows=4,columns=4,#1]
		\begin{puzzleforeground}
			\framepuzzle
			\framearea{black}{(1,1)--(3,1)--(3,3)--(1,3)--cycle}
			\framearea{black}{(3,3)--(3,5)--(5,5)--(5,3)--cycle}
		\end{puzzleforeground}
	}{\end{logicpuzzle}}

\definecolor{ResultColor}{gray}{0.9} 
\definecolor{Dad}{gray}{0.96}  
\definecolor{Medha}{gray}{1.0} 

\newcommand{\bbm}{\begin{boxedminipage}{6.41in}}
	\newcommand{\ebm}{\end{boxedminipage}}

\usepackage{xparse}
\usepackage{lipsum}

\def\exampletext{Example} 
\def\probtext{Problem} 

\NewDocumentEnvironment{testexample}{ O{} }
{
	\colorlet{colexam}{red!55!black} 
	\newtcolorbox[use counter=testexample]{testexamplebox}{%
		empty,
		title={\exampletext\ \thetcbcounter: #1},
		attach boxed title to top left,
		minipage boxed title,
		boxed title style={empty,size=minimal,toprule=0pt,top=4pt,left=3mm,overlay={}},
		coltitle=colexam,fonttitle=\bfseries,
		before=\par\medskip\noindent,parbox=false,boxsep=0pt,left=3mm,right=0mm,top=2pt,breakable,pad at break=0mm,
		before upper=\csname @totalleftmargin\endcsname0pt, 
		overlay unbroken={\draw[colexam,line width=.5pt] ([xshift=-0pt]title.north west) -- ([xshift=-0pt]frame.south west); },
		overlay first={\draw[colexam,line width=.5pt] ([xshift=-0pt]title.north west) -- ([xshift=-0pt]frame.south west); },
		overlay middle={\draw[colexam,line width=.5pt] ([xshift=-0pt]frame.north west) -- ([xshift=-0pt]frame.south west); },
		overlay last={\draw[colexam,line width=.5pt] ([xshift=-0pt]frame.north west) -- ([xshift=-0pt]frame.south west); },%
	}
	\begin{testexamplebox}}
	{\end{testexamplebox}\endlist}

\NewDocumentEnvironment{probexample}{ O{} }
{
	\colorlet{colexam}{red!55!black} 
	\newtcolorbox[use counter=probexample]{probexamplebox}{%
		empty,
		title={\probtext\ \thetcbcounter: #1},
		attach boxed title to top left,
		minipage boxed title,
		boxed title style={empty,size=minimal,toprule=0pt,top=4pt,left=3mm,overlay={}},
		coltitle=colexam,fonttitle=\bfseries,
		before=\par\medskip\noindent,parbox=false,boxsep=0pt,left=3mm,right=0mm,top=2pt,breakable,pad at break=0mm,
		before upper=\csname @totalleftmargin\endcsname0pt, 
		overlay unbroken={\draw[colexam,line width=.5pt] ([xshift=-0pt]title.north west) -- ([xshift=-0pt]frame.south west); },
		overlay first={\draw[colexam,line width=.5pt] ([xshift=-0pt]title.north west) -- ([xshift=-0pt]frame.south west); },
		overlay middle={\draw[colexam,line width=.5pt] ([xshift=-0pt]frame.north west) -- ([xshift=-0pt]frame.south west); },
		overlay last={\draw[colexam,line width=.5pt] ([xshift=-0pt]frame.north west) -- ([xshift=-0pt]frame.south west); },%
	}
	\begin{probexamplebox}}
	{\end{probexamplebox}\endlist}

\newcommand{\bprob} [1]{	\begin{probexample} #1 \end{probexample} }

\usepackage{xspace}
\usepackage[export]{adjustbox}
\newcommand{\Poincare}{Poincar\'e\xspace}
\newcommand{\Dwarf}[2]{\includegraphics[width=#1\textwidth]{#2}}

\begin{document}


\title{On teaching mathematics to gifted students:  some enrichment ideas
	and educational activities}

\author{
\name{Alok Shukla}
\affil{Mathematical and Physical Sciences division, \\School of Arts and Sciences, \\Ahmedabad University, India}
}

\maketitle

\begin{abstract}
Many mathematicians find mathematics aesthetically beautiful and even comparable to art forms such as music or painting. On the other hand, every year a great number of school students leave mathematics with total disillusionment and bitterness, without ever witnessing any beauty in it. In this work, we give some strategies to teach mathematics, especially to gifted students and to instill a love of mathematics in them. We will describe an integrated approach to teaching mathematics where students are introduced to not only more advanced, but also more elegant and more beautiful aspects of the subject early. The proposed integrated approach takes advantage of the fascinating interconnections between various subfields of mathematics, and even borrows from seemingly advanced topics such as number theory and topology. Combining the teaching of computer programming with the teaching of mathematics is another
key focus of this work. This opens up the door to explore, not only the beautiful topics such as fractals, computer math art and computer graphics in general, but also real-life applications resulting from computer simulations in engineering and physics and other natural sciences. We will also discuss the use of storytelling, explorations and experimenting, puzzles and creative problem solving
etc., to make learning interesting for students. It is important to show students that the `true' essence of mathematics goes beyond the dry procedural drill of learning arithmetic. Students must be provided with the opportunities to experience the `aha moment' resulting from the joy of solving a difficult problem, or from the understanding of a deep concept with complete clarity. We will discuss several
strategies, and provide many examples to illustrate how this could be achieved. Some of these strategies also draw from recent advancements in the field of cognitive neuroscience. They also take into account the emotional and psychological aspects of mathematical cognition and learning.
\end{abstract}


\section{Introduction.}
\label{intro}
Imagine the most beautiful, enchanting and colorful garden on the earth ever. Suppose this garden is located on a mountain top, and it's full of different kinds of exquisite flowers, all of vibrant colors. Amongst the flowers some exotic birds are chirping and flying around singing gleefully soft melodious songs, while a few quiet ones are sitting completely still on the branches of nearby lush green trees, as if contemplating the deepest secrets of the universe, and some others are fluttering and dancing carelessly in the sprinkling water coming out of a water-fountain in the middle of the garden. 

Who would not want to visit such a magical garden? However, there is a catch. Before entering the  garden, a visitor has to undertake a very long and arduous journey to the mountain top.
Moreover, visitors must remain blindfolded throughout the climb,  not even once, the visitors are allowed to have a peek into the delightful beauty of the garden. Most of the people could not endure the hardships of the journey and abandoned it in the middle, forgoing their desire to ever visit the garden. They return home, disillusioned and bitter. It is heard in the city that the garden, although created with the utmost architectural precision and perfection, actually looks ugly, and moreover, it has no real-life use for a person living in the city. 

In our opinion, the above analogy sums up one of the most common problems, that many students face in learning mathematics. In many parts of the world, for the first few years in schools, the focus of mathematics education is entirely on teaching arithmetic. Students learn procedural drills of arithmetic, such as how to perform a long division. No doubt, numeracy or number skills are important and all the basic arithmetic operations must be mastered by students. However, a chief problem here is - the almost exclusive focus on arithmetic, and too often on just procedures and drills without a true understanding of the underlying concepts. Unfortunately, to a great extent traditional mathematics teaching  consists of stand-alone skills development through drills and practice exercises. Many students, specially gifted ones,  get bored by the mindless and excessive drills of arithmetic they are often subjected to. It is like blindfolding them and not allowing them to sneak a peek into the beautiful garden of mathematics. After several years of such education, it is hardly surprising if only a few students still remain interested in learning mathematics.

In this work, we will focus on ways to address this issue and propose some strategies to teach mathematics to gifted students.
First, we recall some traditional learning models. According to Swiss psychologist Jean Piaget, cognitive developments in children take place in four distinct stages: sensory-motor, preoperational, concrete, and formal \cite{piaget1964part}. A progression through these stages is supposed to be age-dependent and linear. Therefore, the model posits that a child of a certain age can only deal with concrete physical objects and is yet unable to conceptualize abstractly. 

Here we differ a little from Piaget's learning models while addressing the problem of devising mathematics course plans for gifted students. In our opinion, gifted children should not be forced to learn mathematics linearly. Mathematics is not a straight ladder that children should climb step by step. Indeed, mathematics is full of beautiful interconnections. Children should not be prevented from learning the so-called `difficult' topics beyond their ages. We have found impressive learning capabilities and imaginations in young children, and if nurtured appropriately they have the potential to learn advanced mathematical concepts.

As we have remarked earlier, spending too much time on procedural aspects of arithmetic alienates many students. Such students perceive mathematics as a boring subject and lose their motivation to learn it. Therefore, all efforts should be made to introduce beautiful parts of mathematics to children, as early as possible. Of course, in following this non-linear approach to learning, there might remain gaps in their understanding of concepts, and they need to go back and forth to learn some of the missing concepts. Still, this way of learning is preferable in many cases and it offers several advantages. Various subfields in Mathematics are beautifully interconnected, and sometimes these interconnections help in understanding a mathematical concept by providing views from multiple angles. Therefore, an early introduction to algebra, geometry, and coordinate geometry is very helpful to gifted students in their mathematical journey.

\begin{tcolorbox}
	We strongly recommend an integrated approach of teaching mathematics with an \emph{early} introduction to: 
	\begin{itemize}
		\item  algebra, abstraction and symbolism in mathematics,
		\item  geometry and coordinate geometry,
		\item  combinatorics and graph theory,
		\item  beauty in mathematics,
		\item  computer programming.
	\end{itemize}	
	The proposed integrated approach will also take into account:
	\begin{itemize}
		\item  puzzles and problem-solving strategies, 
		\item  independent mathematical exploration,
		\item  applications in physics and engineering
		\item  managing emotional well-being and dealing with failures.
	\end{itemize}	
\end{tcolorbox}

In the following, we will describe our key ideas in detail and illustrate them with examples, with topics ranging from preschool to high school mathematics, and sometimes even beyond that. We will refer to research in cognitive neuroscience, and discuss, how both the educators and the learners can make use of results in this field. Above everything else, developing a deep love for mathematics in young minds should be the ultimate goal of mathematics education. Once this goal is achieved, the other pieces of the puzzle will automatically fall in place.

This work contains several examples, including some python codes. We note that many of the concepts discussed here were tested on a set of gifted students. The names of these students appear in the examples to follow. We are thankful to these ideal and enthusiastic students. Teaching these students provided opportunities for the author to learn and reflect on mathematics and its teaching. The ages of these students in the context of examples in this work can be treated as Natasha - 4 years, Sophie - 6 years, Anya and Medha - 7 years, and Pragya and Arya - 10 years. Natasha has yet to start her formal schooling. Sophie was being home-schooled. Everyone else was a member of the GT (Gifted \& Talented) program in a public school in Oklahoma, United States. \\

\noindent \textbf{A Remark.} Even though, we have used the term  `gifted students'  throughout this work, we do not want to imply that this work is meant exclusively for only a very few students. Without engaging in any debate over `nature' vs `nurture', we would like to state that, in our opinion, most of what we have to say in this work applies to any school student of mathematics, who is willing to devote the time and efforts needed in the pursuit of learning this beautiful subject.

\section{Let us count.}  \label{Sect:Counting}
Learning to count is amongst the first mathematical skills that a child learns. Counting is an important survival skill. Even in primitive societies, infants and non-human animals, counting skills are developed \cite{beran1998chimpanzee}, \cite{boysen2014counting}, \cite{davis1982counting}. How does a child learn to count? Research suggests that some aspects of counting abilities are innate \cite{antell1983perception}. However, in a formal learning environment most often a child learns to count by actually associating `words' to the count of discrete objects such as toys, cups or fruits. 

What could be a good strategy to teach counting and basic addition facts to young children? Before attempting to answer this question, we recall that the human brain is biologically so well-wired for image processing as against number crunching. In numerical calculations, such as multiplying or dividing two large numbers, or finding the cube root of a large prime number,  even a calculator can easily beat us. However, until very recently humans were doing better than the computer in image recognition and analysis \cite{fleuret2011comparing}, although in near future perhaps computers will catch up and even surpass humans in image processing as well. \footnote{In fact, already in 2015 researchers at Baidu Research, claimed that on one of the computer vision benchmarks, the ImageNet classification challenge, their system has exceeded the human recognition performance \cite{wu2015deep}.} Still, a human brain has this special ability to complete even a partially hidden picture and attach meaning to it. Perhaps once upon a time survival of humans depended upon quickly identifying a predator, even with incomplete visual data, and promptly acting. As a result, even a five old kid is capable of identifying and placing the missing piece of a puzzle and thus completing the picture. Well, if we humans are so well-equipped to carry out visual image processing, then why not try to include visual aids wherever possible in education, including in learning how to count and how to learn basic addition facts?

\begin{figure}
	\begin{subfigure}{.45\textwidth}
		\centering  	
		\includegraphics[scale=0.25, frame]{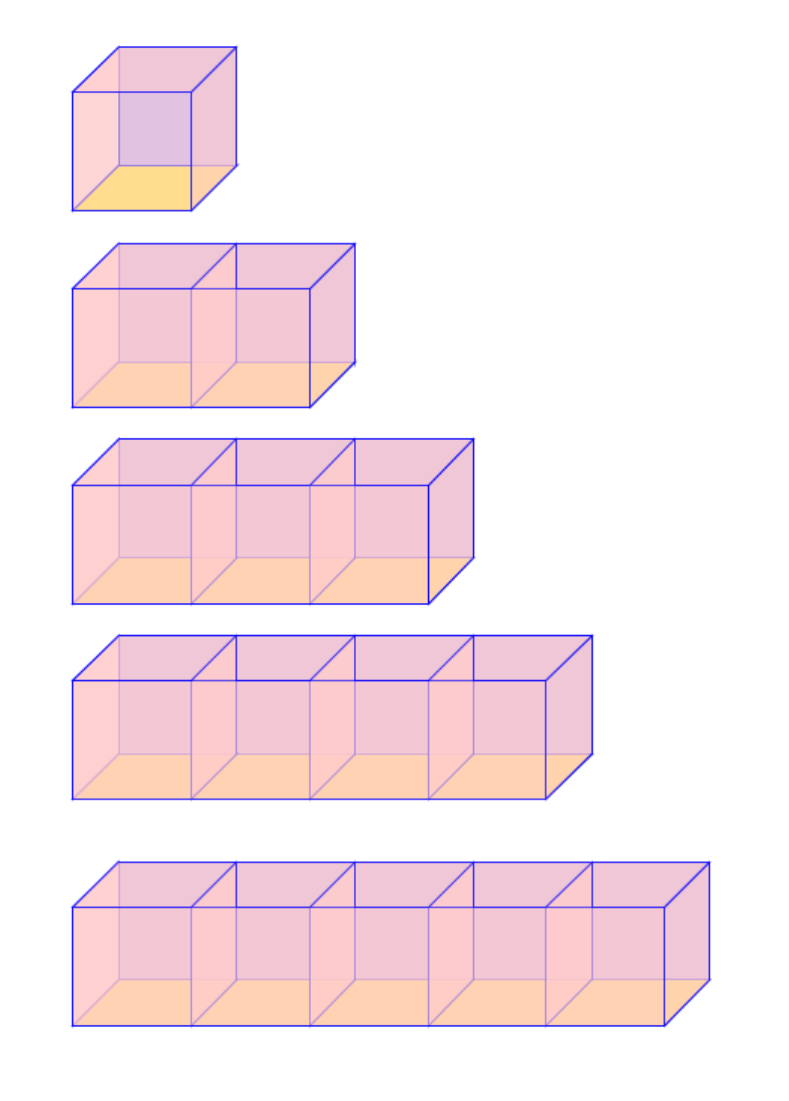}
		\caption{} \label{fig:blocks}
	\end{subfigure}
	\hspace{1 cm}
	\begin{subfigure}{.45\textwidth}
		\centering
		\includegraphics[scale=0.25,frame]{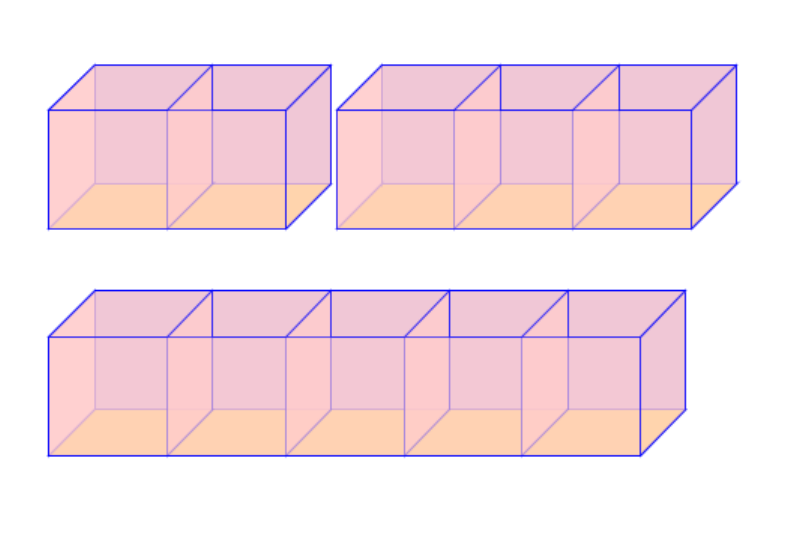}
		\caption{} \label{fig:blocksadd}
	\end{subfigure}
	\caption{Developing an 'intuitive' number-sense in students using visual aids. } 
\end{figure}

Keeping the above discussion in mind, in addition to discrete counting objects, students should also be taught counting using continuous models using detachable cubic blocks, as shown in Fig.~$ \ref{fig:blocks}$. Any other `continuous' objects or household items may also be used. In our experience with $ 3$-$4 $ years old children, after a few trials, they could easily visually identify the blocks by just looking at their sizes. For example, Natasha could identify the blocks after playing and practicing with them for a few days.  Later, when she was asked to identify the block $ 3 $, she immediately responded three, without counting. In fact, she could identify any individual block of sizes $ 1 $ to $ 5 $ without counting. The biggest advantage of this visual analog approach was in learning addition facts. As Natasha could identify the blocks of sizes $ 2 $ and $ 3 $, and also she knew that when she joined them, then the new joined block was as big as a block of size $ 5 $. Therefore, she instinctively knew that $ 2 + 3 =5 $ (see Fig.~$ \ref{fig:blocksadd} $).  

\begin{figure}
	\mbs
	\begin{subfigure}{.45\textwidth}
		\centering  	
		\includegraphics[scale=0.4]{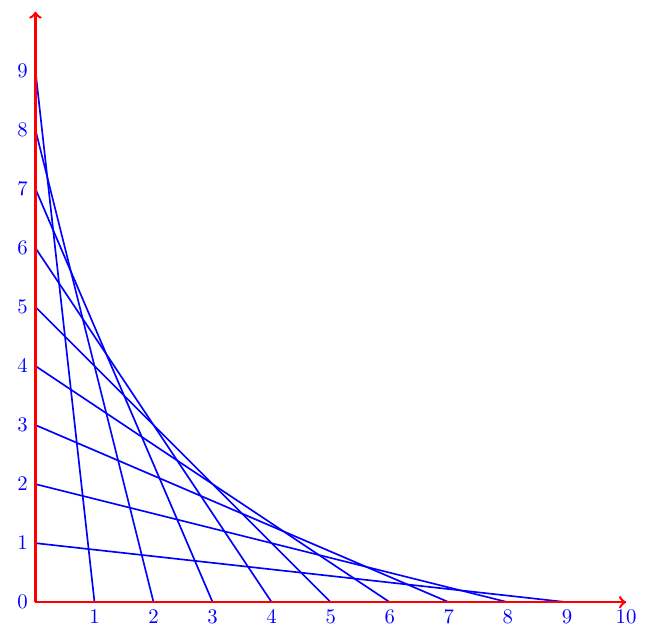}
		\caption{} \label{fig:sum10}
	\end{subfigure}
	\hspace{1 cm}
	\begin{subfigure}{.45\textwidth}
		\centering
		\includegraphics[scale=0.3]{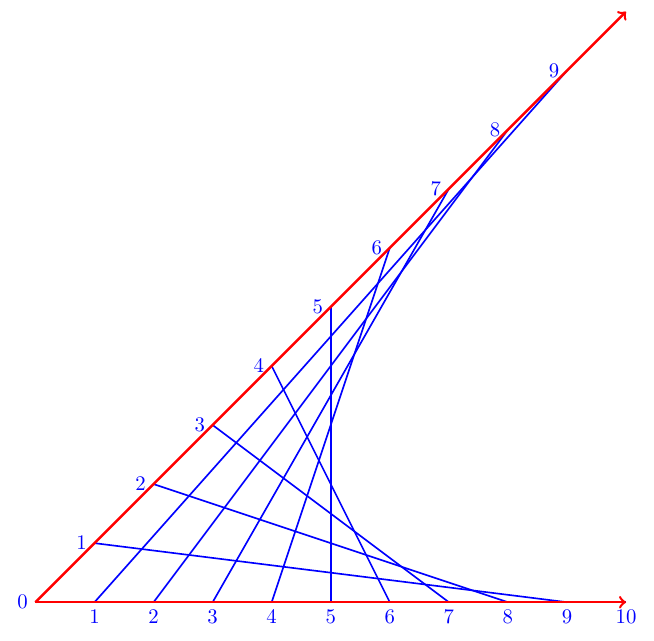}
		\caption{} \label{fig:sum10xy}
	\end{subfigure}
	\vspace{1cm}
	\begin{subfigure}{.45\textwidth}
		\centering  	
		\includegraphics[scale=0.22]{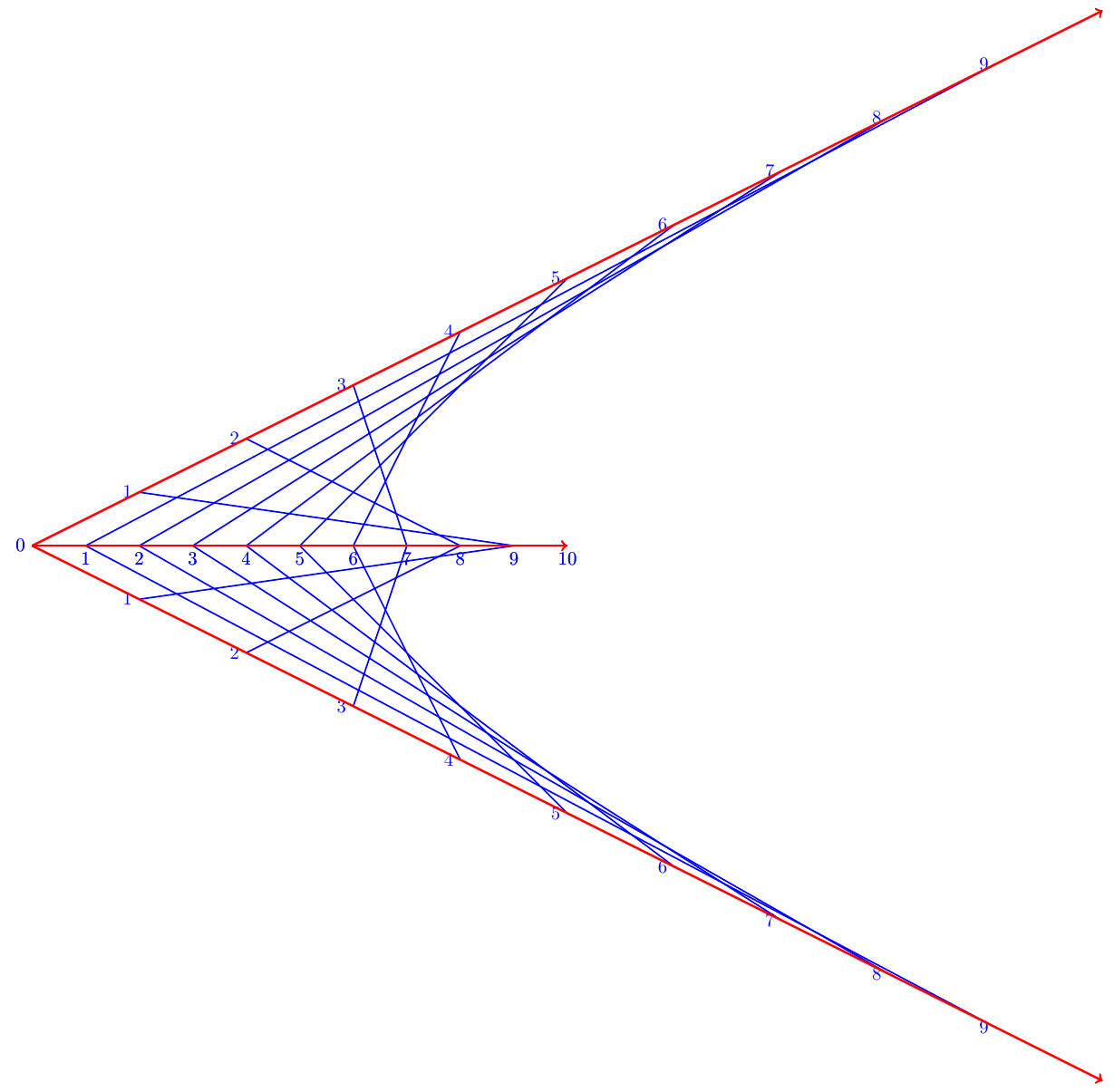}
		\caption{} \label{fig:sum10xy2}
	\end{subfigure}
	\hspace{1 cm}
	\begin{subfigure}{.45\textwidth}
		\centering
		\includegraphics[scale=0.015]{S100xy2.pdf}
		\caption{} \label{fig:sum100xy2}
	\end{subfigure}
	\mbe	
	\caption{Curve stitching.}\label{fig:sumexamples}
	
\end{figure}

Later in this work, we will discuss the importance of engaging students using mathematical games, beautiful patterns and hands-on activities. Here we present one example of a beautiful pattern that Natasha obtained by joining points on two perpendicular lines such that their sum was $ 10 $ (see Fig.~$ \ref{fig:sum10} $).   She then carried out some variations of this exercise to draw some other interesting patterns (see Fig.~$ \ref{fig:sum10xy} $, \ref{fig:sum10xy2} and \ref{fig:sum100xy2}). She first performed these exercises on paper. Later, on she did it on a computer using the interactive software Geogebra (see \url{https://www.geogebra.org/}). We remark here that, such hands-on activities with visualization help in developing the so-called intuitive `number-sense' in gifted students.

\section{A relationship game: meeting abstractions early.} \label{Sect:Relationship}
Mathematics is concerned with relationships among various objects. Even before a child learns how to count, this aspect could be emphasized. For instance, when asked - ``Who is the father of your dad?'', Natasha immediately responded - ``Grandpa''. She could not answer at first when she was asked a slightly complicated relationship question - ``Who is the wife of your father`s dad?'',  but after a little help, she answered the question correctly. Within a few days of practicing with such kinds of problems, she could not only answer even more complicated relationship problems but also challenged others with her own questions. It became a sort of game that she enjoyed playing.

The relationship game that Natasha played helped set the stage for the next level. Once she was familiar with the relationship game, she was given problems like -``Suppose D is C's daughter. Also, suppose A's father is B, and B's mother is C. A and E are B's daughters. Then who is A's grandmother?''. 


\begin{figure}
	\centering
	\includegraphics[scale=0.75,frame]{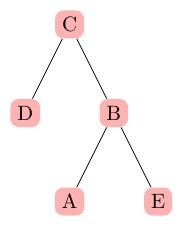}
	\caption{A family tree.} \label{fig:a_family_tree}
\end{figure}

She learned how to draw a family-tree (see Fig.~$ \ref{fig:a_family_tree}$) and then she was able to correctly answer such questions. She understood that we can represent people by their names or a symbol, which is just a shorter name, for convenience.  She also learned how to represent the relationships amongst people, by connecting the symbols using lines with arrow heads. She was already getting introduced to some abstractions.   

Another occasion to learn abstract symbolism came when Natasha was learning how to add numbers. She was given a series of problems, and she learned the addition facts such as 
$ 5 +0 = 5$, $6 +0 =6, \ldots 100 + 0 =100  $ and so on. At the end of a series of such problems, she was asked -``What would be $ x + 0 $?'', she responded within a few seconds $ x $. She didn't even ask - ``What does $ x $ stand for?''.  Maybe the relationship game she had played earlier, helped her. When probed further to explain why should $ x + 0  $ be $ x $, she did not clearly state that there the $x $ stood for any arbitrary number, instead she responded -``0 is nothing, so adding it to something should work like that''. Clearly, she could notice the pattern and could come to a conclusion that demanded mathematical abstraction and generalization. Later on, when Natasha learned subtraction, she could correctly answer that $ x- x $ should be equal to $ 0 $. Similarly, when Natasha was learning multiplication, after answering a series of questions related to mathematical facts  such as $ 1 \times 2 =2 $, $ 1 \times 3 =3 $, and $ 1 \times 4 =4 $; she could also answer the question $ 1 \times x = x $. 

The above examples show that at least some kids will be able to understand and appreciate mathematical abstractions from an early age. We should not wait until students gain  complete mastery of arithmetical concepts before exposing them to higher-level abstractions and algebraic thinking. In fact, there are ways to integrate arithmetic teaching with algebraic thinking at each stage. We will discuss some other examples later in this work.
\section{Beauty in mathematics.} \label{Sect:Beauty}
\begin{figure}
	\mbs
	\begin{subfigure}{.45\textwidth}
		\centering  	
		\includegraphics[scale=0.4]{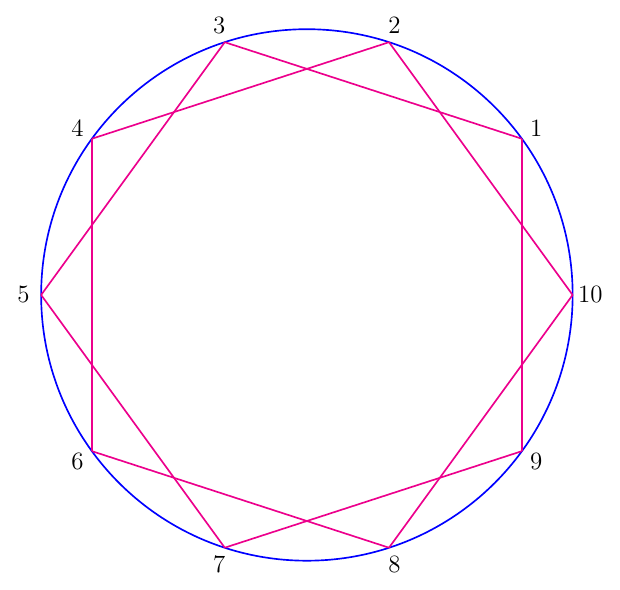}
		\caption{} \label{fig:Skip1}
	\end{subfigure}
	\hspace{1 cm}
	\begin{subfigure}{.45\textwidth}
		\centering
		\includegraphics[scale=0.4]{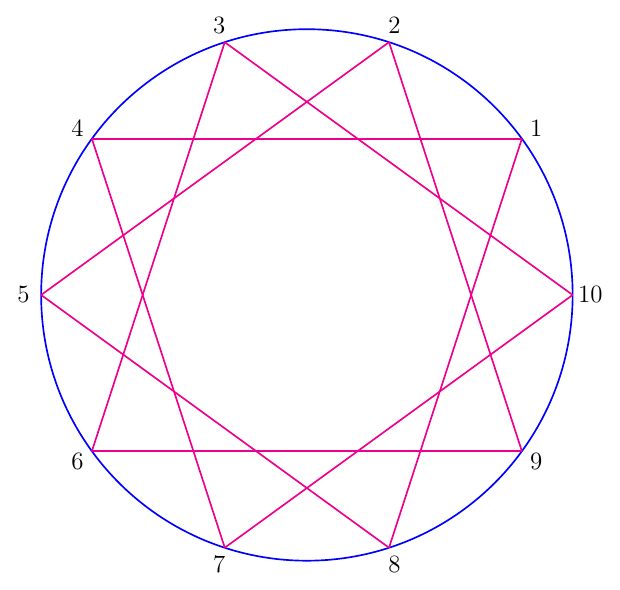}
		\caption{} \label{fig:Skip2}
	\end{subfigure}
	\vspace{1cm}
	\begin{subfigure}{.45\textwidth}
		\centering  	
		\includegraphics[scale=0.4]{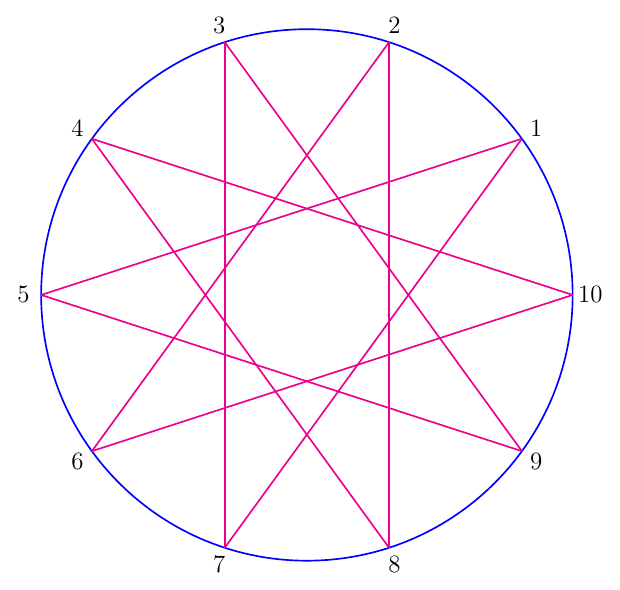}
		\caption{} \label{fig:Skip3}
	\end{subfigure}
	\hspace{1 cm}
	\begin{subfigure}{.45\textwidth}
		\centering
		\includegraphics[scale=0.4]{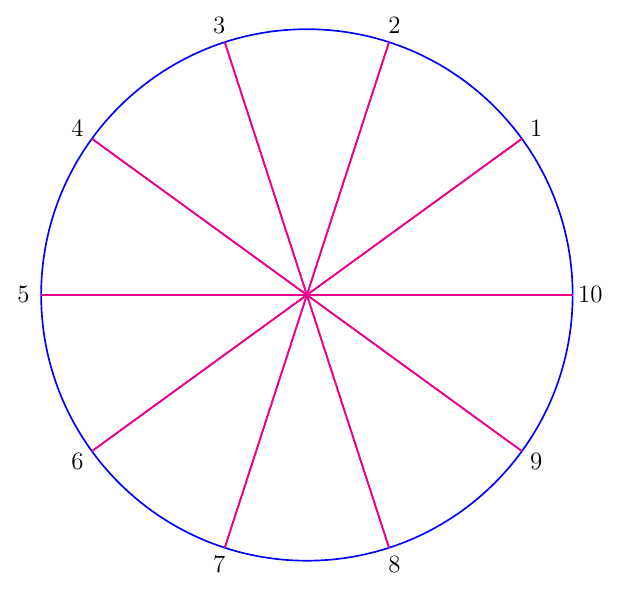}
		\caption{} \label{fig:Skip4}
	\end{subfigure}
	\mbe
	\caption{Skip counting on a circle.}\label{fig:skipexamples}
\end{figure}

Many accomplished mathematicians have commented on the beauty of mathematics. The following quote by \Poincare  succinctly describes the importance of beauty in mathematics from the perspective of a pure mathematician: ``The mathematician does not study pure mathematics because it is useful; he studies it because he delights in it, and he delights in it, because it is beautiful.''  
\begin{tcolorbox}
	A recent research result in neuroscience, \cite{zeki2014experience}, shows that the experience of mathematical beauty and the experience of visual, musical, and moral beauty are correlated by the activity in the same brain area, the medial orbito-frontal cortex.	
\end{tcolorbox}
It is important to introduce young students to the aesthetic beauty of mathematics as soon as possible. Even while students are learning basic arithmetic operations such as addition, skip counting, multiplications etc., it is possible to include beautiful visual patterns or math art in their learning modules. We have already discussed one such example in Sect.$ \ref{Sect:Counting} $ (see Fig.~$ \ref{fig:sumexamples}$). Here, we present some more examples. Medha enjoyed the following activity. Suppose, the numbers $ 1$, $ 2,  \ldots, 10 $ were marked on a circle such that the numbers are equally spaced (see Fig.~$ \ref{fig:skipexamples} $). Now suppose a turtle starts from the number $ 1 $ and visits numbers on the circle in the counterclockwise fashion by skipping the next number in the sequence. Clearly, such a turtle will travel along the path $ 1-3-5-7-9-1 $. Similarly, if the turtle starts from $ 2 $ then it will travel along the path $ 2-4-6-8-10-2 $ (see Fig.~$ \ref{fig:Skip1}$). Therefore, Fig.~$ \ref{fig:Skip1}$ is the pattern obtained by performing a `skip-counting' by one on a circle with ten numbers. Similarly, Medha obtained Fig.~$ \ref{fig:Skip2}$, Fig.~$ \ref{fig:Skip3}$ and Fig.~$ \ref{fig:Skip4}$ by performing skip counting by $ 2,3 $ and $ 4 $ respectively.

We note that Fractals are another example of beautiful mathematical objects, that are even capable of imitating the shapes of many natural objects.  Sophie was intrigued by the beauty of fractal patterns when she first encountered them. She learned how to draw 
a Sierpinski triangle (Fig.~$ \ref{fig:Sierpinski} $) and Koch snowflakes.

The advantages of introducing such beautiful objects to young children are manifold.  First, children are naturally curious and once they find fractals beautiful they want to discover more about fractals and get interested in mathematics. Second, the colorful fractal patterns present the opportunity for teachers to engage students by letting them create beautiful fractal arts by coloring and possibly cutting shapes out of paper. Third, for advanced students fractals could be used to teach various mathematical concepts. Fractals can be used to familiarize students with the basics of complex numbers and the concept of recursion. After gaining some knowledge of computer programming, students can write computer programs to generate fractal patterns. In Sect.~$ \ref{Sect:Fractals} $ and Sect.~$\ref{Sect:programming_physics}$, we will further describe using fractals in the context of computer programming.

\section{Math and play.} \label{Sect:MathandPlay}
Children are playful by nature. Play is important for the cognitive, social, emotional and physical well-being of children \cite{milteer2012importance}. Playing provides a great opportunity for parents and teachers to connect to a child. Many mathematical concepts can be introduced to children by playing with them. Pragya was hopping over alternate squares while saying the numbers $ 2$, $ 4 $, $6,\ldots$, aloud. She repeated this activity with $3$, $ 6 $, $9,\ldots  $, and learned the tables of $ 2 $ and $ 3 $ while hopping over the squares. Pragya also liked to play the board game `Snake and Ladders' with her sister. She also played a modified version of this game with a dice with numbers $ 10$, $ 20 $, $30,\ldots,60 $ on it. Playing such games helped her in mastering basic counting skills. 

A $ 4 \times 4 $ Sudoku can be turned into a game for young children, wherein players take turn to solve the puzzle by placing the fruits on the Sudoku board placed on the floor.

\subsection{Nim}
Combinatorial games provide a rich source of puzzles for developing problem-solving skills in young children. Suppose two players play the following game. They take turns to eat $ a = 5 $ apples and $ b = 6 $ bananas. On her turn, a player has to select a type of fruit that she wants to eat. Then she must eat at least one fruit of that type. If she wants, she may eat more than one fruit of the chosen type. The player who eats the last remaining fruit wins the game (See Table~$ \ref{Tab_Nim} $). This is a two `heap' version of the famous combinatorial game Nim (See \cite{bouton1901nim}). The good thing about this game is that it has been completely solved, in the sense that there exists an algorithm for giving the winning strategy for one of the players from any specified game position.  For the $ a \times b $ version of the game described earlier, finding this winning strategy is not so difficult, and discovering this winning strategy could provide a good intellectual challenge to a gifted child. 

\newcommand{\apple}[1]{
	\foreach \i in {1,...,#1}{
		\tikz{\includegraphics[scale=0.05]{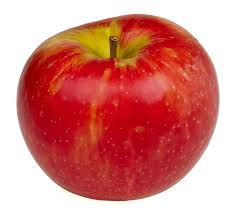}} \hspace{0.2cm}  }}

\newcommand{\banana}[1]{
	\foreach \i in {1,...,#1}{
		\tikz{\includegraphics[scale=0.05]{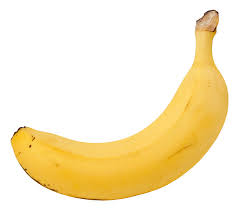}} \hspace{0.2cm}  }}

\begin{table}[hbt]
	\centering
	\label{tbl:heatwheel_res}
	\begin{tabular}{lll}
		\toprule
		\textbf{Move} &
		\textbf{Apples} &
		\textbf{Bananas}\\
		\midrule 
		Initial position \qquad & \apple{5} \qquad & \banana{6}     \\ 
		Player $ 1 $ eats $ 3 $ bananas \qquad & \apple{5} \qquad & \banana{3}  \\ 
		Player $ 2 $ eats $ 2 $ apples \qquad & \apple{3} \qquad & \banana{3}  \\ 
		Player $ 1 $ eats $ 1 $ apple \qquad & \apple{2} \qquad & \banana{3}  \\ 
		Player $ 2 $ eats $ 2 $ bananas \qquad & \apple{2} \qquad & \banana{1}  \\ 
		Player $ 1 $ eats $ 1 $ apple \qquad & \apple{1} \qquad & \banana{1}  \\ 
		Player $ 2 $ eats $ 1 $ apple \qquad &  \qquad & \banana{1}  \\ 
		Player $ 1 $ eats $ 1 $ banana \& wins the game \qquad &  \qquad &   \\ 
		\bottomrule
	\end{tabular}
	\caption{5 apples and 6 bananas game.} \label{Tab_Nim}
\end{table}

For example, Pragya enjoyed playing this game. After a while, she started to notice some of the attributes of the winning position. With a little help from her teacher, she listed down the winners for the small games. For the $ (a=1) \times (b=1) $ game, it was clear to her that the second player will be the winner, irrespective of what the first player did. Next, for the $ (a=2) \times (b=1)  $ game, the first player can just eat one apple to convert this game to the $ (a=1) \times (b=1) $ game, with the first player now effectively becoming the `second player', and thereby winning the game as discussed earlier. Pragya noticed that the same strategy worked for $ (a >1) \times (b=1) $ games, wherein the first player eats $ a-1 $ apples to convert the game to $ (a=1) \times (b=1) $ game. What about a $ (a=2) \times (b=2) $ game? It was clear that because of the symmetry one needed to analyze the only two moves available to the first player. The first, eating the $ 2 $ apples, which is a losing move as then the second player eats $ 2 $ bananas. The second, eating $ 1 $ apple converts the game to the $(a=1) \times (b=1)   $ game, which was analyzed earlier. But this was also a losing move for the first player. Therefore, the first player will lose the $ (a=2) \times (b=2) $ game if the second player plays perfectly. Now Pragya could see a pattern emerging. It was clear to her that for the $ (a=3) \times (b=2) $ game, the first player will simply eat one apple and reduce the game to the $ (a=2) \times (b=2) $ game with the first player now becoming the second. It was obvious that the first player had a winning strategy for the $ (a=3) \times (b=2) $ game.  At this point, Pragya conjectured that if $ a = b $ then the second player will have a winning strategy, otherwise the first player has a winning strategy. In fact, after a little more exploration and analysis she discovered that if $ a=b $, then the goal of the second player is to always keep this balance on her turn, i.e., a winning strategy for the second player is to follow the moves of the first player but for the other fruit. In the case of $ a \neq b $, without loss of generality suppose $ a >b $. Then the first player will eat $ a-b $ apples, reducing it to the $ b \times b $ game, becoming effectively the second player, and thereby winning it. After understanding the winning strategy, Pragya went on to play this game with her younger sister Medha. She was happy to win several games, until Medha quit playing. In fact, Medha was also noticing some patterns then, and she went to her teacher for discussing the game with him.  She was given the following problem to think over.

\bprob{Suppose two players take turns to eat $ a = 5 $ apples, $ b = 6 $ bananas and $ c = 7 $ oranges. On her turn a player must select and eat one type of fruit. She is free to eat as many available fruits of the chosen type  as she wants. The player who eats the last remaining fruit wins the game. Is there a winning strategy for any player in this game? If yes, describe this winning strategy?}

\subsection{Make $ 10 $}
For students who are learning basic addition facts, the following game could be useful. The aim is to create a pile of $ 10 $ apples. Two players take turns placing either one apple or two apples in the pile. The player who makes the last move wins the game. A sample game is shown in Table $ \ref{Tab:make10} $.

\bprob{Does there exist a winning strategy for any player in the ``Make $ 10 $'' game? If yes, what is it?}

Of course, there can be many variations of this game, such as ``Make $ 15 $'' or ``Make $ 100 $''. Also, instead of only one or two apples, players might be allowed to place one, two or three apples, or any other such variations.

\bprob{Does there exist a winning strategy for any player in the ``Make $ 15 $'' game? If yes, what is it?}

\begin{table}[hbt]
	\centering
	\label{tbl:heatwheel_res}
	\begin{tabular}{ll}
		\toprule
		\textbf{Move} &
		\textbf{Apples} \\
		\midrule 
		Initial position \qquad &     \\ 
		Player $ 1 $ puts $ 1 $ apple in the pile \qquad & \apple{1} \\ 
		Player $ 2 $ puts $ 2 $ apples in the pile \qquad & \apple{3} \\ 
		Player $ 1 $ puts $ 2 $ apples in the pile \qquad & \apple{5} \\ 
		Player $ 2 $ puts $ 1 $ apple in the pile \qquad & \apple{6} \\ 
		Player $ 1 $ puts $ 1 $ apple in the pile \qquad & \apple{7} \\ 
		Player $ 2 $ puts $ 1 $ apple in the pile \qquad & \apple{8} \\
		Player $ 1 $ puts $ 2 $ apples in the pile \& wins! \qquad & \apple{10} \\ 
		\bottomrule
	\end{tabular} 
	\caption{Make $ 10 $ game. Two players take turns to place either one apple or two apples in the pile. The player who plays the last move wins the game.} \label{Tab:make10}
\end{table}

\section{Combinatorics and graph theory.}\label{Sect:Combinatorics}
Combinatorics and graph theory are two important subjects that are not included in a typical middle school curriculum. We believe that these two subjects should be taught to  gifted students in elementary and middle schools. The fundamental principle of counting is pretty intuitive and can be visually explained to a child who knows how to multiply. For example,  Sophie could understand and correctly answer the following problem. 
\say{There are 3 flights from California to France and 2 flights from France to India. Sophie wants to fly from California to France and then to India. How many different ways are there for her to fly from California to India?}. 

\begin{figure}
		\centering
		\includegraphics[scale=0.6] {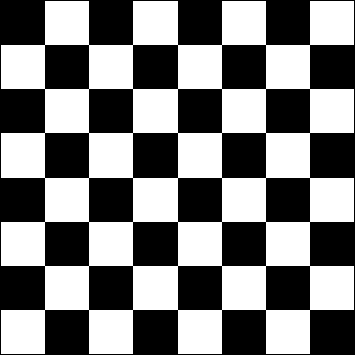}
	\caption{In how many ways can eight rooks be placed on a chessboard so that no two rooks are attacking each other?} \label{fig_chess_rooks}
\end{figure}

The subject of combinatorics offers many attractive problems which could provide intellectual nourishment to a gifted child. For example, the following problem stumped Pragya, and she spent quite some time trying to solve it. In how many ways can eight rooks be placed on a chessboard (Fig.~$ \ref{fig_chess_rooks} $) so that no two rooks are attacking each other?

Many problems in Graph theory can be stated and understood by even elementary school kids, however, they are often hard to prove. For example, students could be encouraged to explore and conjecture / discover the \say{four color theorem} (see \cite{gonthier2008formal}). Another attractive problem to challenge young children is the ``K\"{o}nigsberg bridge problem'' (also see Sec.~$ \ref{Sect:Puzzle} $). Graph theory can be used as a vehicle to introduce the concept of \say{proof} in mathematics to students. Pragya was introduced to some basic definitions in Graph theory and then given a set of theorems and results to prove. 
\begin{itemize}
	\item {Prove that the sum of the degrees of the vertices of any finite graph is even.}
	\item Prove that a finite simple graph, with more than one vertex, has at least two vertices with the same degree.
	\item Prove that a finite graph has an even number of vertices with odd degrees.
	\item Show that a connected graph on $ n $ vertices is a tree if and only if it has $ n - 1 $ edges.
	\item Prove that a complete graph on $ n $ vertices contains $ \displaystyle {\frac{n(n-1)}{2}} $ edges.
	\item Prove that for a given vertex of an odd degree in a graph, there exists a path from it to another odd-degree vertex in the graph.
\end{itemize}

Some of the problems were too difficult for Pragya. However, with some helpful hints, she could solve the above problems. She was asked to write her solutions, as a professional mathematician would write. Unfortunately, writing and communicating mathematics clearly, is a skill that is often not emphasized in schools.
\subsection{Adjacency Matrix and counting the number of walks.}
\begin{figure}
	\centering
	\includegraphics[scale=0.75,frame]{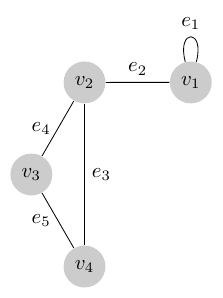}
	\caption{Counting the number of walks of a given length.} \label{fig_walk_adjacency}
\end{figure}

Matrices are one of the most important concepts in mathematics. Of course, linear transformations and vector spaces are the right settings to introduce them to a gifted student. However, since in this section we are discussing graphs, we would like to point out that young students can be taught about computing the number of walks of a given length from the adjacency matrix. Pragya and Medha both learned how to multiply two matrices. Then they were taught to write the adjacency matrix for a given graph. For a given graph with its vertex set $ V = {v_1,v_2,\cdots,v_n}$, its adjacency matrix $ A $ is a square matrix of size $ n $ such that the $ (i,j)^{\text{th}} $ entry is given by the number of edges directly connecting the vertex $ v_i $ to $ v_j $.  Loops are also counted once, i.e., for a loop at $ v_i $ (when $ v_i $ is connected to itself by an edge), the  $ (i,i)^{\text{th}} $ diagonal entry will be $ 1 $.  Pragya and Medha both easily learned how to calculate the adjacency matrix of a given graph. For example, for the graph in Fig.~$\ref{fig_walk_adjacency}  $ they calculated the adjacency matrix $ A $ and its powers $ A^2 $, $ A^3 $ and $ A^4 $ as follows
\begin{align*}
&A = \left[\begin{smallmatrix}
1 & 1 & 0 & 0 \\
1 & 0 & 1 & 1 \\
0 & 1 & 0 & 1 \\
0 & 1 & 1 & 0
\end{smallmatrix}\right], \quad A^2 = \left[\begin{smallmatrix}
2 & 1 & 1 & 1 \\
1 & 3 & 1 & 1 \\
1 & 1 & 2 & 1 \\
1 & 1 & 1 & 2
\end{smallmatrix}\right], \quad  A^3 = \left[\begin{smallmatrix}
3 & 4 & 2 & 2 \\
4 & 3 & 4 & 4 \\
2 & 4 & 2 & 3 \\
2 & 4 & 3 & 2
\end{smallmatrix}\right],
 \,\text{and} \quad 
A^4 =
\left[\begin{smallmatrix}
7 & 7 & 6 & 6 \\
7 & 12 & 7 & 7 \\
6 & 7 & 7 & 6 \\
6 & 7 & 6 & 7
\end{smallmatrix}\right].
\end{align*}
Now they were asked to count the number of  walks from $ v_2 $ to $ v_2 $ of  length $ 2 $ (in a walk a vertex or an edge can be repeated and the length of a walk is the number of edges traversed counting with its multiplicity). They listed the $ 3 $ such possible walks: $ e_2e_2 $, $ e_3e_3 $ and $ e_4e_4 $. Of course, the $ (2,2)^{\text{th}} $ entry of $ A^2 =3 $. Similarly, they listed the number of  walks from $ v_2 $ to $ v_2 $ of  length $ 4 $ and found it amazing that the number of such walks is $ 12 $ which is the  $ (2,2)^{\text{th}} $ entry of $ A^4$. They had fun listing other walks and how their numbers are given by the entries of the matrices listed above. They could conjecture that the number of possible walks of length from $ v_i $ to $ v_j $ is given by $ (i,j)^{\text{th}} $ entry of $ A^n $. Of course, it was a great puzzle for them to figure out why is it so!

\section{Distributive law.} \label{Sect:Distributive}
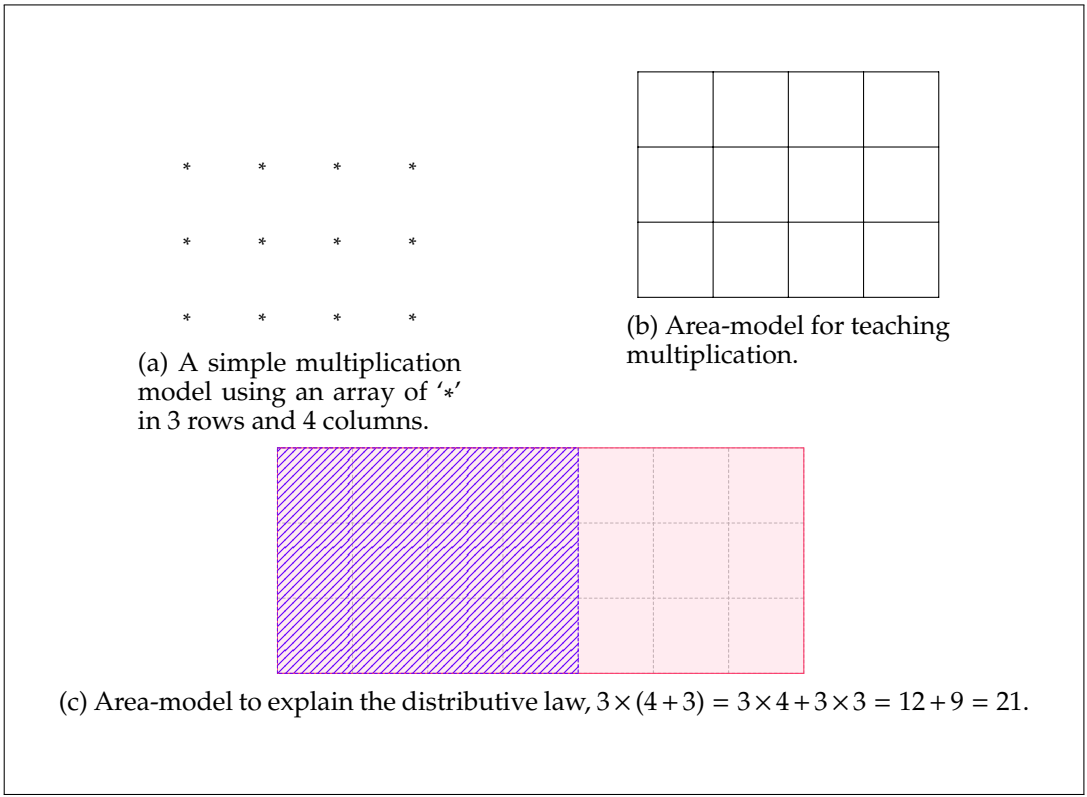
\begin{figure}
	\centering
	\begin{subfigure}{.3\textwidth}
		\centering
		\includegraphics[scale=0.6]{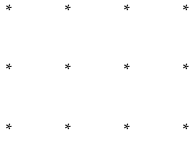}
		\caption{{\footnotesize A simple multiplication model using an array of `$ * $' in $ 3 $ rows and $ 4 $ columns.}}	\label{fig_star_model}
	\end{subfigure} \hspace{2 cm}
	\begin{subfigure}{.3\textwidth}
		\centering
		\begin{tikzpicture}[scale=0.6,line cap=round,line join=round,>=triangle 45]
		\draw [ xstep=1.0cm,ystep=1.0cm] (0,0) grid (4,3);
		\clip(0,0) rectangle (4,3);
		\end{tikzpicture}
		\caption{{\footnotesize Area model for teaching multiplication.}}\label{fig_area_model_first}	
	\end{subfigure}
	\vspace{0.5 cm}
	\begin{subfigure}{.9\textwidth}
		\centering  	
		\includegraphics[scale=0.6]{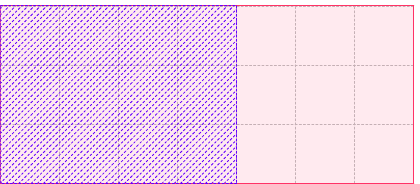}
		\caption{{\footnotesize Area model to explain the distributive law, $ 3 \times (4 +3) = 3 \times 4 + 3 \times 3  = 12 + 9 = 21$.}}\label{fig_area_model_distributive}	
	\end{subfigure}
	\caption{Distributive law and the area model.} 
\end{figure}

It is very convenient to have real-world models for mathematical operations such as multiplication. These models can then be used to introduce further mathematical symbolism and abstraction. A simple multiplication model is given in Fig.~$ \ref{fig_star_model} $. Another simple multiplication model is the area model given in Fig.~$ \ref{fig_area_model_first} $. 
Traditionally, two and higher digit multiplication is taught by using the column method. The traditional method is a fast and efficient method of doing multiplication by hand. However, often students do not really understand why the multiplication algorithm works. They just follow it mechanically without any actual insight into the whole process. On the other hand, using the distributive property to teach multiplication to young children has the advantage that the concepts are understandable with an easy visualization offered by the `area model' of the multiplication operation (Fig.~$ \ref{fig_area_model_distributive} $).
Moreover, on following this approach students learn the concepts of `area' and `multiplication' together, unlike the traditional method wherein the concept of area is introduced much later. Interestingly, the distributive property illustrates an important problem-solving strategy in mathematics (as well as maybe in our lives): when one is unable to solve a hard problem, one should try to break it into smaller problems that are easier to solve. This concept can be used to teach how to perform multiplication as well as division operations. Anya initially had some trouble learning multiplication and distributive law using the area model. However, once she understood this method, she preferred it over the column method of multiplication, that she was learning at her school.

\begin{figure}
	\centering
	\includegraphics[scale=0.6]{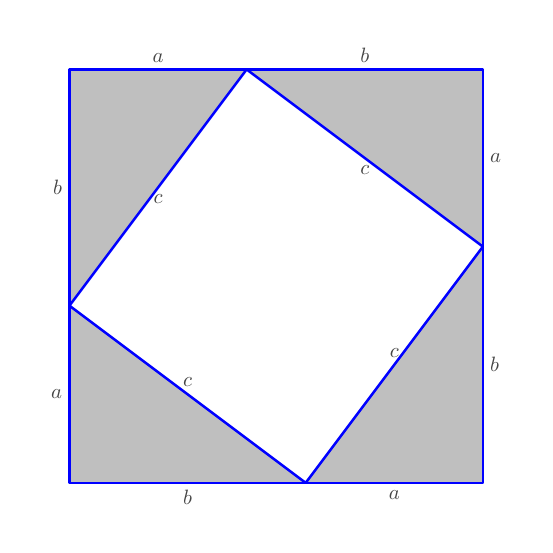}
	\caption{Garfield's proof of the Pythagorean theorem.} \label{fig:Pyhtagorus_proof}
\end{figure}

Another advantage of introducing the distributive law  and the area model is that it helps students in understanding the distributive law $ a(b+c) = ab + ac $, when they encounter it in algebra later. In fact, the Pythagorean theorem, arguably one of the most important theorems in mathematics, can be proved using the distributive law and the area model. Such a proof of the Pythagorean theorem was given by the twentieth president of the United States James A.~Garfield. Medha could easily understand a modernized version of Garfield's proof. Consider the big square in Fig.~$ \ref{fig:Pyhtagorus_proof} $. We compute its area in two different ways. First, its area is $ (a+b)^2 $. Then its area is also the sum of the areas of four triangles which equals $ 4 \times \frac{1}{2} a b = 2 ab$ and the area of a small square $ c^2 $. Therefore, we must have
\begin{align*}
(a+b)^2 = 2 ab +c^2 \implies a^2 + b^2 = c^2. 
\end{align*}  

\section{Early introduction to algebraic symbols and functional thinking.} \label{Sect:EarlyAlg}
Algebra is one of the biggest hurdles that many college students in the United States face during their college mathematics coursework. We refer readers to \cite{stewart2017algebra} for an interesting take on the consequences of algebra under-performances at the college level. We believe that one solution to this problem is an early introduction to algebraic thinking in school education. Indeed, algebraic symbols can be introduced even before preschool as we have remarked earlier. The same is true for the concept of functions. Functions are arguably one of the most important concepts in entire mathematics. Many five-year-olds are capable of understanding the fundamental idea behind functions. Natasha was given a paper with the names of her sisters and friends on the one side and their ages on the other (see Fig.~$ \ref{Figure_functionn_ages} $), and she was asked to draw arrows connecting names to corresponding ages. She could draw it correctly. Later on, she was introduced to the function notation and she understood the meaning of $ f(\text{Pragya}) = 11 $. Of course, she also understood that since two persons can have the same age, it is okay if two or more arrows from the \say{left} point to the same object on the \say{right}. However, since one person can`t have two different ages at the same time, from one object on the \say{left}
more than one arrow can't originate. Later on, she was given several other examples of functions, and she even came up with a few examples of her own. 

One particular example, related to the \say{relationship games} she played earlier (see Sect.~\ref{Sect:Relationship}) pertained to the composition of functions. Suppose $ f $ is a map from the son to his father, i.e., $ f(A) = B $ if $ A $'s father is $ B $. Then, Natasha understood that the statement \say{A`s father is B, and B`s father is C} can be expressed in the function-notation as $ f(f(A)) =C $.

Another example simply related to the function $ f(x) =x+1 $. She was given the information that $ f(1)=2,f(2)=3,f(3)=4 $ and so on. Then she was asked to guess the pattern and predict $ f(9) $, she could get the correct answer $ 10 $. However, when she was asked what should be $ f(x) $ she said $ y $. She reasoned that $ y $ comes after $ x $, so $ f(x) $ should be $ y $. Perhaps this misconception originated from her misunderstanding of the \say{ordinal} and \say{cardinal} numbers. Moreover, she was also unaware of the symbolic expression \say
{$ x+1 $}. She was then introduced to another game. Some marbles were hidden under a cup  and $ 1  $ marble was placed beside it on the table. She was told that there are $ x $ marbles hidden inside the cup and if she adds the $ 1 $ marble beside it the total would be $ x+1 $ marbles. Then the cup was lifted to show that $ x $ was indeed $ 2 $, i.e., $ 2 $ marbles were hidden. And, in that case there were total $ x+1 = 2+1 =3 $ marbles. She also played another variation of this game. This time there were two cups, labeled $ x $ and $ y $ hiding some unknown number of marbles. Her goal was to count the total number of hidden marbles, i.e. $ x+y $. She turned the cup labeled $ x $ and noted the number of hidden marbles under the heading $ x $. Similarly, she counted the number of marbles hidden under the cup labeled $ y $ and noted it under the heading $ y $. Finally, she computed $ x + y $. After several variations of this game, she understood what was going on. 
The upshot of all the above exercises was that not only she became familiar with the function notation and simple symbolic expressions, but she also learned how to model and express her own thoughts using the function notation. This early introduction to functional thinking will be helpful to her in her mathematical journey in the coming years.

\begin{figure}
	\centering
	\begin{tikzpicture}[scale =0.4,
	>=stealth,
	bullet/.style={
		fill=black,
		circle,
		minimum width=1pt,
		inner sep=1pt
	},
	projection/.style={
		->,
		thick,
		shorten <=10pt,
		shorten >=10pt
	},
	every fit/.style={
		ellipse,
		draw,
		inner sep=0.8pt,
		fill = purple,
		opacity = 0.3
	}
	]
	\foreach \y/\l in {1/Natasha,2/Medha/,3/Pragya,4/Sophie,5/Anya}
	\node[bullet,label=left:$\l$] (a\y) at (0,2*\y) {};
	
	\foreach \y/\l in {1/4,2/8,3/11,4/6,5/10}
	\node[bullet,label=right:$\l$] (b\y) at (15,2*\y) {};
	
	\node[draw,fit=(a1) (a2) (a3) (a4) (a5),minimum width=4.3cm] {} ;
	\node[draw,fit=(b1) (b2) (b3) (b4) (b5),minimum width=3cm] {} ;
	
	\draw[projection] (a1) -- (b1);
	\draw[projection] (a2) -- (b2);
	\draw[projection] (a3) -- (b3);
	\draw[projection] (a4) -- (b4);
	\draw[projection] (a5) -- (b2);
	\end{tikzpicture}
	\caption{The function $ f $ mapping names of Natasha and her friends to their ages. } \label{Figure_functionn_ages}
\end{figure}

\section{Modular arithmetic.} \label{Sect:Modular}
Gauss has called number theory the `queen of mathematics'. Indeed, number theory is a source of fascinating problems, which can often be formulated in a relatively simple language, for example, the famous modularity theorem of Wiles, i.e., Fermat's last theorem. The modular arithmetic developed by Gauss is a beautiful topic in number theory, which is easily accessible to a gifted student in  middle school. It is also amazing that modular arithmetic has recently found applications in creating error correcting codes, for example, Reed-Solomon error correction codes, which is useful in modern devices such as  DVDs, Blu-ray Discs and  even in satellite communications. This is just one example where pure mathematics pursued only for the sake of its beauty and human curiosity turns out to be useful in practical applications.

The example of a clock is a well-known device to introduce modular arithmetic to beginners. Pragya very easily understood the concept of addition modulo 12. It was clear to her that for a $ 12  $-hour clock 
$ 1 $ is the same as $ 13 $.

\begin{figure} 
	\centering
	\includegraphics[scale=0.6]{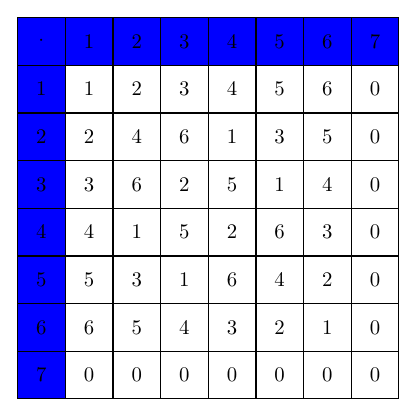}
	\caption{Multiplication table modulo $ 7 $.} \label{fig:multiplication_modulo_seven}
\end{figure}

Then for simplicity, she was given a few numbers and asked to divide them by $ 3 $ and find remainders if any. After this exercise, she observed that when a number is divided by $ 3 $, the only possible remainders are $ 0 $, $ 1 $ or $ 2 $. In this manner, the set of integers, which is an infinite set, is partitioned into $ 3 $ classes (equivalent classes), with a number $ n $ being in the class corresponding to the remainder obtained when $ n $ is divided by $ 3 $. Next Pragya learned how to add or multiply two numbers modulo $ 3 $ and how these operations are independent of representatives of the three classes chosen. She then created a multiplication table modulo $ 3 $ and modulo $ 7 $ (see Fig.~$ \ref{fig:multiplication_modulo_seven} $). Arithmetic Modulo $ 7 $ is connected with the days of the week, just as arithmetic modulo $ 12 $ is connected with a clock, and so it constitutes another familiar example of modular arithmetic. In fact, she was also introduced to the definition of a Group and asked to check if the additions and multiplications modulo $ 7 $ satisfy group axioms. Moreover, using the multiplication table modulo $ 7 $ she could discover interesting number theoretic facts such as $ 6^{2n} -1 $ is divisible by $ 7 $ for each natural number $ n $.

\begin{center}
	\begin{figure} 
		\mbs
		\begin{subfigure}{.5\textwidth}
			\centering  	
			\includegraphics[scale=0.5]{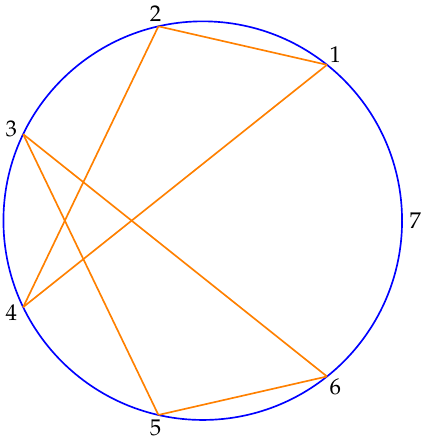}
			\caption{}\label{fig:sub-first}
		\end{subfigure}
		\begin{subfigure}{.5\textwidth}
			\centering  
			\includegraphics[scale=0.5]{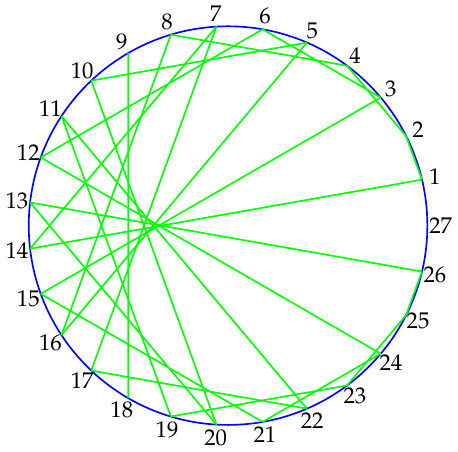}
			\caption{} 	\label{fig:sub-second}
		\end{subfigure}
		\vspace{1 cm}
		\begin{subfigure}{.5\textwidth}
			\centering  	
			\includegraphics[scale=0.5]{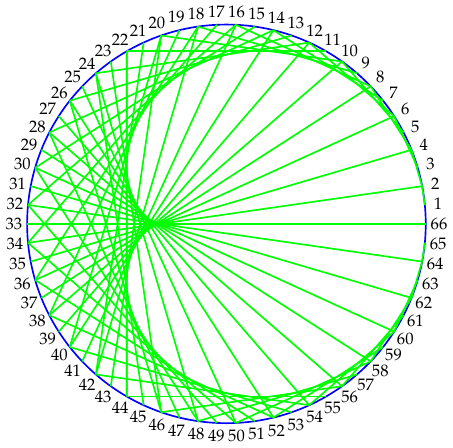}
			\caption{} \label{fig:sub-third}
		\end{subfigure}
		\begin{subfigure}{.5\textwidth}
			\centering  	
			\includegraphics[scale=0.5]{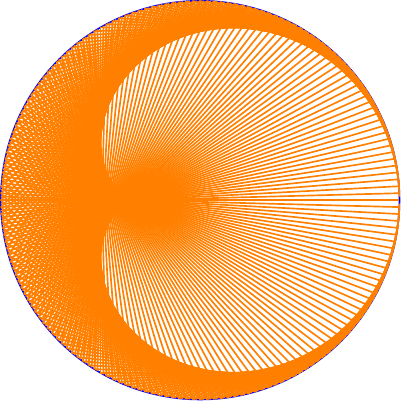}
			\caption{} \label{fig:sub-fourth}
		\end{subfigure}	
		\mbe
		\caption{Modular arithmetic and cardioid. Points $ 1, \, 2, \ldots n$ are placed on a circle with equal spacing between points. Then chords are drawn between points $ m $ and $ (2m \mod n) $ for each point $ m $ on the circle. The diagrams (a), (b), (c) and (d) are drawn using $ n=7$, $ 36 $, $66 $ and $ 360 $ respectively. }
	\end{figure}
\end{center}

Another interesting exercise for Pragya was to use the above multiplication table to create a digram as follows. She started with drawing a circle. Then she placed the points $ 1$, $ 2, \ldots, n=7 $ on the circle with equal spacing between points. Now for each of the above points $ m $, she joined the points $ m $ and $ (2m \mod 7) $ by a straight line to obtain Fig.~\ref{fig:sub-first}. The interesting thing happened when she did the same exercise for the set $ {1,\, 2,\ldots, n=36} $ to get Fig.~\ref{fig:sub-second}. The envelope of the chords happened to be a `heart like shape'. Indeed, this shape was a cardioid and the above method to draw a cardioid as the envelope of a pencil of lines is due to L.~Cremona. Then it was a good time to introduce polar coordinates to Pragya. She understood how to plot a cardioid $ r =  1 + cos(\theta)$. It also became obvious to Pragya that drawing the chords and doing modular arithmetic is going to be very time-consuming if $ n $ is big. Therefore, it would be nice if we can take help from a computer. This was a good motivation for her to learn to write a computer program (also see Sect.~\ref{Sect:Scratch}). She learned the concept of 'loop' and managed to write programs for generating Figs.~\ref{fig:sub-third} and \ref{fig:sub-fourth}. 

\begin{figure}
	\mbs
	\begin{subfigure}{.5\textwidth}
		\centering  	
		\includegraphics[scale=0.75]{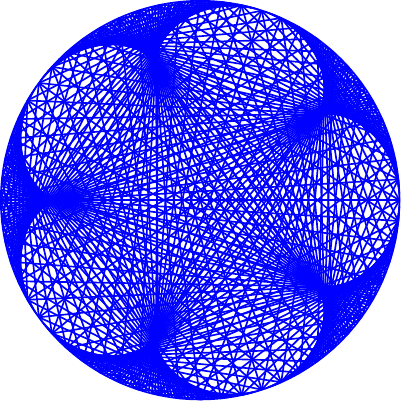}
		\caption{} \label{fig:sub-one-three-sixty}
	\end{subfigure}
	\begin{subfigure}{.5\textwidth}
		\centering  	
		\includegraphics[scale=0.75]{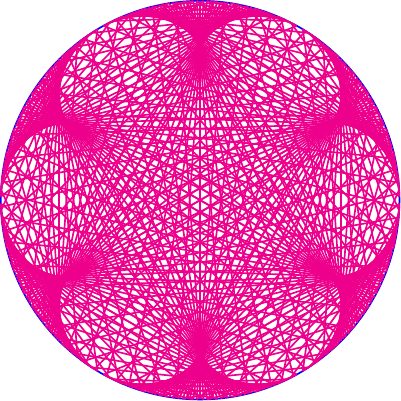}
		\caption{} \label{fig:sub-two-three-sixty}
	\end{subfigure}
	\vspace{1 cm} \\
	\begin{subfigure}{.5\textwidth}
		\centering  	
		\includegraphics[scale=0.75]{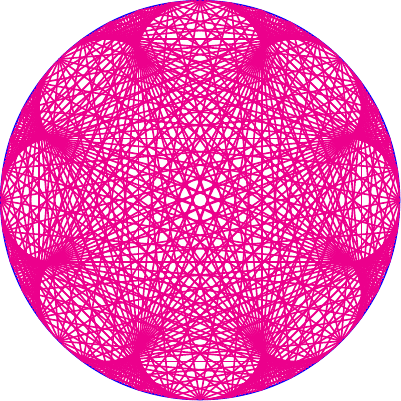}	
		\caption{} \label{fig:sub-three-three-sixty}
	\end{subfigure}
	\begin{subfigure}{.5\textwidth}
		\centering  	
		\includegraphics[scale=0.75]{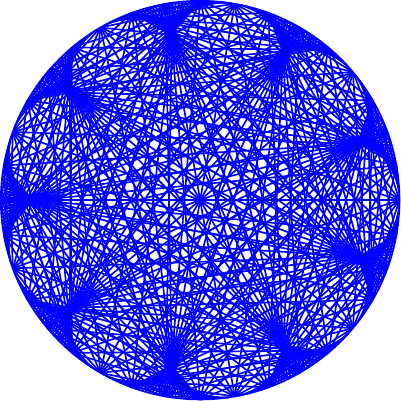}	
		\caption{} \label{fig:sub-four-three-sixty}
	\end{subfigure}
	\mbe
	\caption{Multiplication modulo $ 360 $. The two diagrams on the top row ((a) and (b) are created using multiplication by 6 and 7 modulo 360 respectively. The bottom row diagrams are generated by using multiplication by 9 and 10 modulo 360 respectively. The appearance of geometric shapes such as pentagon, hexagon etc.,  at the centers of the diagrams above is very interesting. } \label{fig:three-thirty-sixty}
\end{figure}

Other patterns can be generated by considering more points on the circle, for example with $ 360  $ points. Patterns generated by considering multiplication modulo $ 360 $ are shown in Fig.~$ \ref{fig:three-thirty-sixty} $. The two diagrams on the top row ((a) and (b) are created using multiplication by $ 6 $ and $ 7 $ modulo $ 360 $ respectively. The bottom row diagrams are generated by using multiplication by $ 9 $ and $ 10 $ modulo $ 360 $ respectively. Interesting geometric shapes such as pentagon, hexagon etc., appear at the centers of the diagrams in Fig.~$ \ref{fig:three-thirty-sixty} $.

\begin{figure}
	\centering
	\includegraphics[scale=0.2]{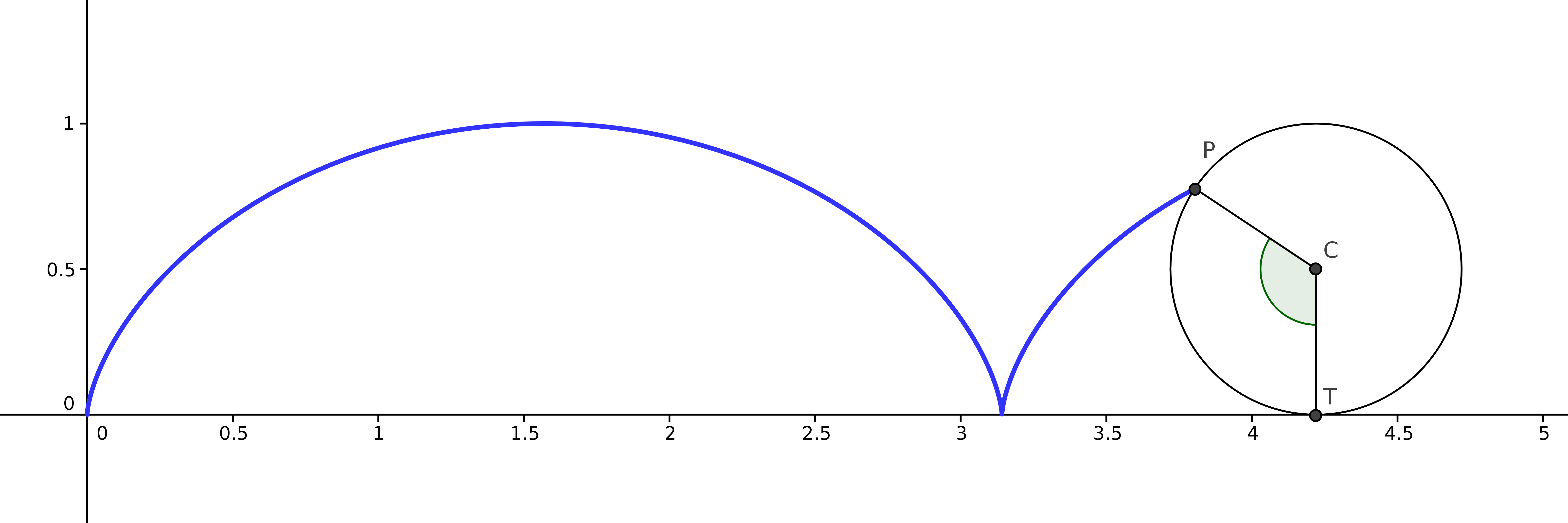}
	\caption{A cycloid. Pick a point P on a  circle, which is rolling over a line without slipping. The locus of the  P is a cycloid. This was simulated using Geogebra. }\label{fig:cycloid}
\end{figure}

Cardioid is a cycloidal curve. Cycloidal curves are described by the path traced by a point that is lying to a circle rolling along another circle or a straight line without slipping (See Fig.~$ \ref{fig:cycloid} $ and Fig.~$ \ref{fig:Cardiod} $). 
Interactive geometry and algebra Software Geogebra, can be used to simulate this rolling and create an animation depicting various cycloidal curve (see Sect.~$ \ref{Sect:Geogebra} $).
Cycloids appear in the solution of the brachistochrone problem (see \cite{de1993galileo} for a historical account of this problem) and Pragya was fascinated to learn the history of the brachistochrone problem and its connection with the cycloids.

\begin{figure}
	\centering
	\includegraphics[scale=0.6]{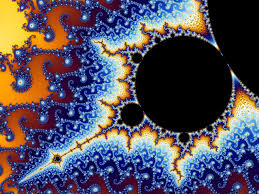}
	\caption{Mandelbrot set. The central bulb has the shape of  a cardioid.} \label{fig:Mandelbrot}
\end{figure}

There was even more to the story of a cardioid.  The complex map $ {\displaystyle z\to z^{2}}  $ sends a circle through the origin to a cardioid. This results in the interesting fact that the  boundary of the central bulb of the Mandelbrot set is a cardioid. The Mandelbrot set contains an infinite number of self-repeating patterns such that the central bulbs of these smaller repeating copies are approximate cardioids (see Fig.~$ \ref{fig:Mandelbrot} $). The Mandelbrot set is an example of Fractals.  Fractals are one of the most visually stunning and mysteriously beautiful objects in mathematics (see Sect.~$ \ref{Sect:Fractals} $).  

\begin{figure}
	\mbs
	\begin{subfigure}{.45\textwidth}
		\centering  	
		\includegraphics[scale=0.06]{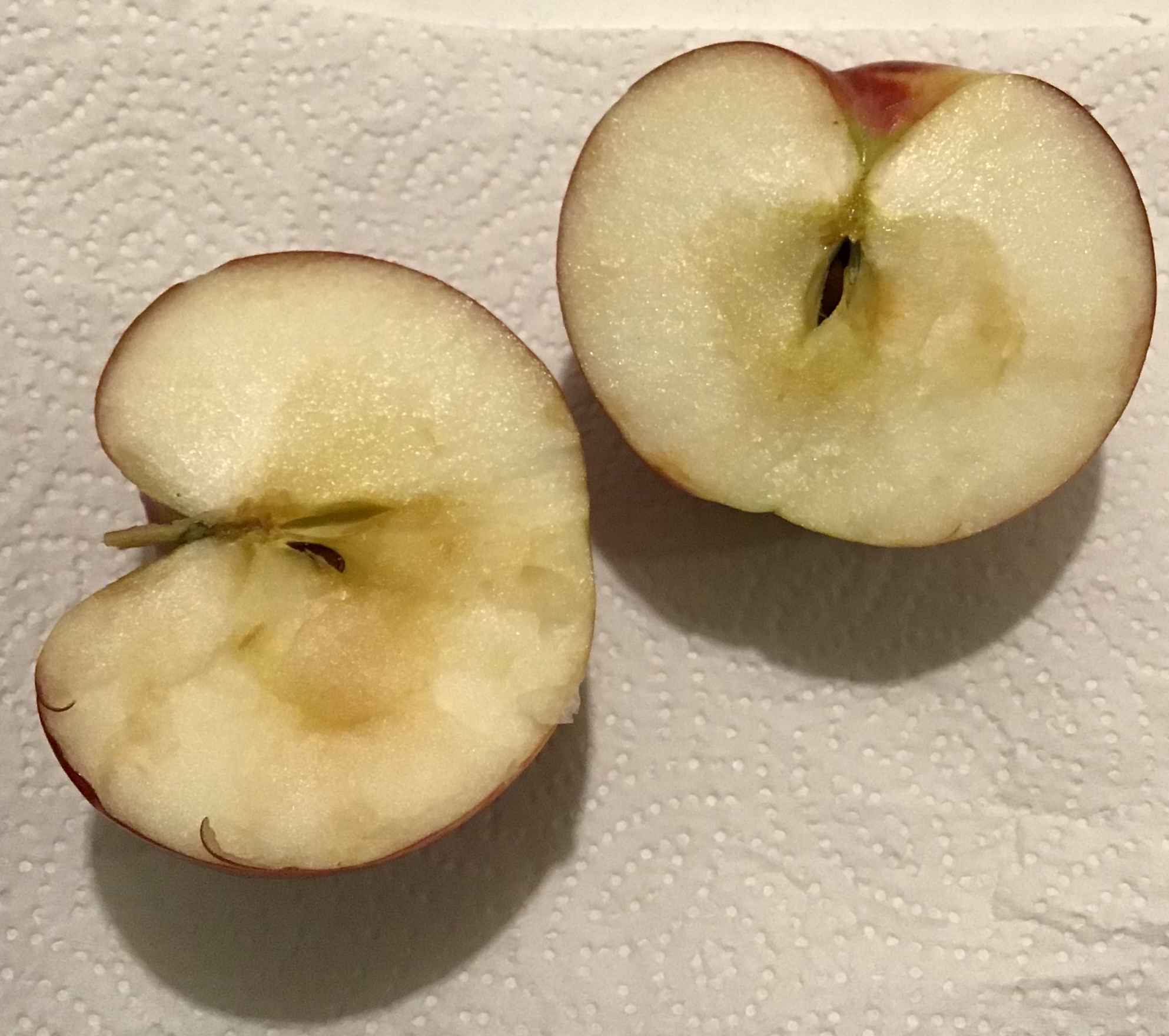}
		\caption{Cardioid and apples.} \label{fig:apple}
	\end{subfigure}
	\hspace{1 cm}
	\begin{subfigure}{.45\textwidth}
		\includegraphics[scale=0.065]{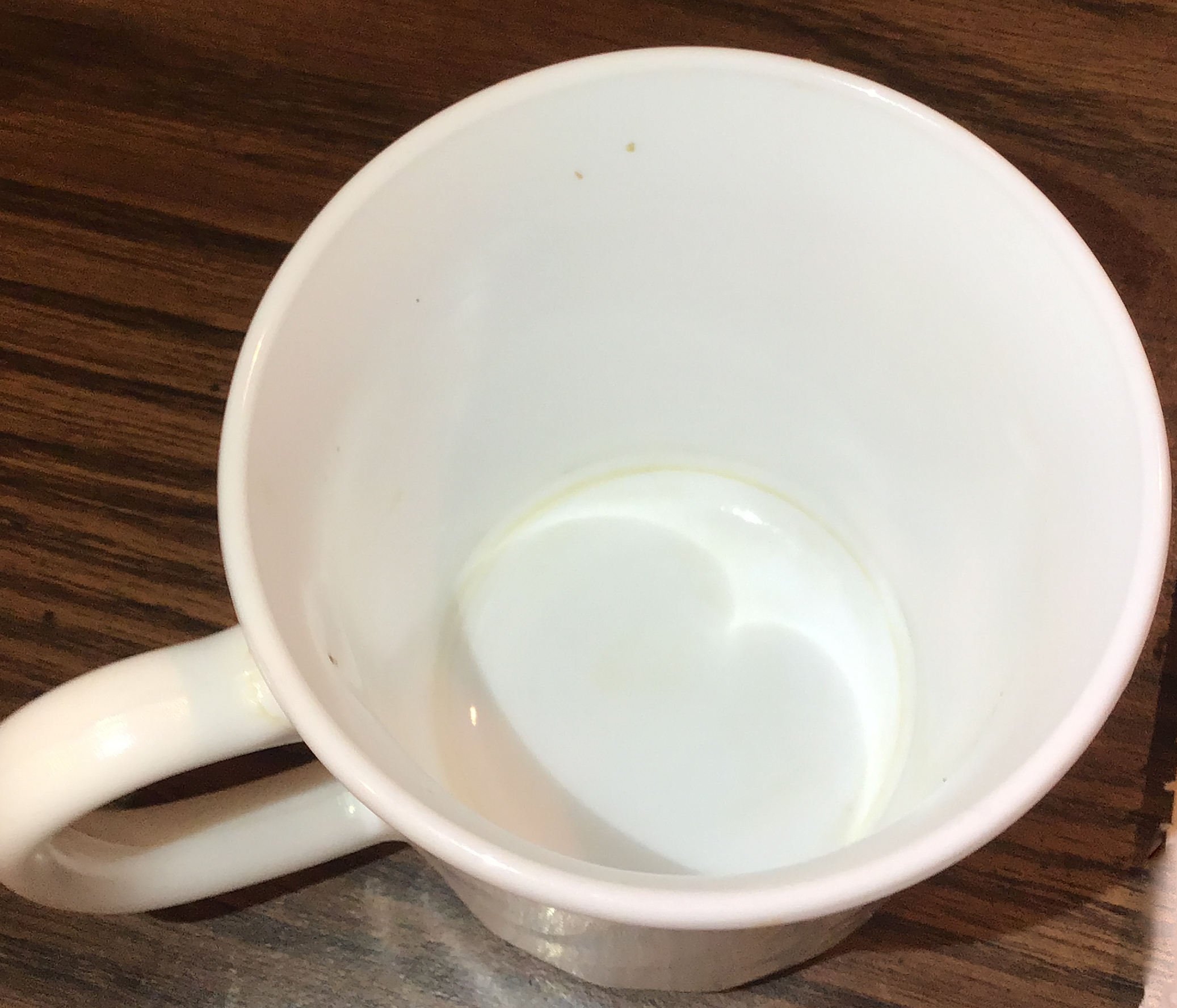}
		\caption{A Cardioid appearing, when a torch was lit in front of a coffee cup.} \label{fig:cofeecup}
	\end{subfigure}
	\mbe
	\caption{Cardioid in real life.}
\end{figure}

A couple of days after learning about the cardioid, Pragya found a piece of one half of an apple and observed that it looked very similar to a cardioid (see Fig.~$ \ref{fig:apple} $).
Later she was amazed to see that when a torch was lit in front of a coffee cup she could see the shape of a cardioid inside it (see Fig.~$\ref{fig:cofeecup}  $). It was a great opportunity for her to learn about the nature of light and how it travels in a medium and the phenomenon of reflection and refraction of light etc.


\section{Storytelling in mathematics.} \label{Sect:Story}
Storytelling is a powerful mode of human communication and it is capable of engaging us at a deep emotional level. Stories often emerge from shared cultural constructs and common human experiences and they foster a unique bond between the storyteller and the listener. Teachers and educators have long advocated in support of the benefits of reading bedtime stories to children. It is believed that storytelling provides many educational benefits to children, such as improvement in vocabulary, and an enhancement in imagination and communication skills. Indeed, stories provide listeners an opportunity to paint the characters on the canvas of their own imaginations and create their own movies in their minds as the story unfolds. Perhaps it explains the findings of a recent research study, \cite{yabe2018effects}, wherein the effects of storytelling on the brains of children were assessed by using near-infrared spectroscopy (NIRS) and it was concluded that -``The results indicated more sustained brain activation to storytelling in comparison with picture-book reading,
suggesting possible advantages of storytelling as a psychological and educational medium in children.'' There are many neuroscience and cognitive science-based research papers that discuss the merits of storytelling as an instructional tool \cite{aldama2015science}. In fact, researchers have studied the use of storytelling, in marketing \cite{pulizzi2012rise},  as a primary leadership tool \cite{7500928}, and as an alternative method of healing for trauma survivors \cite{carey2006expressive}.

\begin{figure} 
	\centering
	\begin{tikzpicture} [scale=0.65]
	\draw[very thick,<-] (-6.3,0)--(6.3,0);
	\draw[very thick,<-] (0,-6.3)--(0,6.3);
	\tkzInit[xmax=6,ymax=6,xmin=-6,ymin=-6]
	\tkzGrid
	\tkzAxeXY
	\node (King) at (0,0) {\Dwarf{0.04}{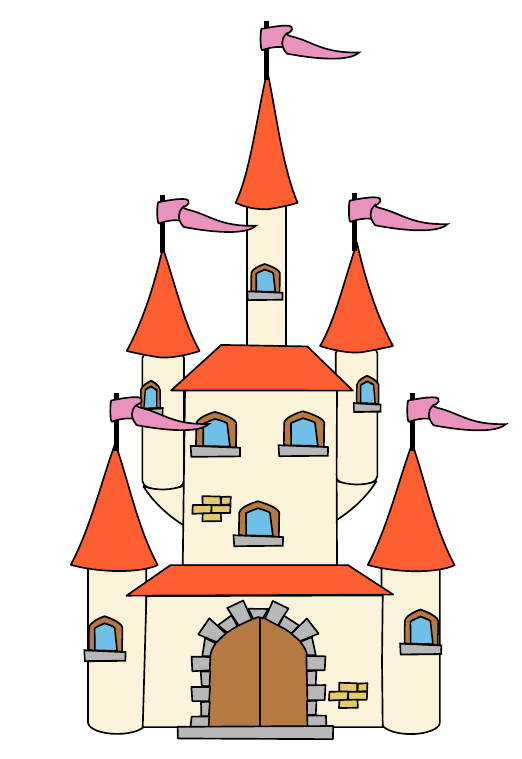}}; 
	\node (Magic) at (3,4) {\Dwarf{0.03}{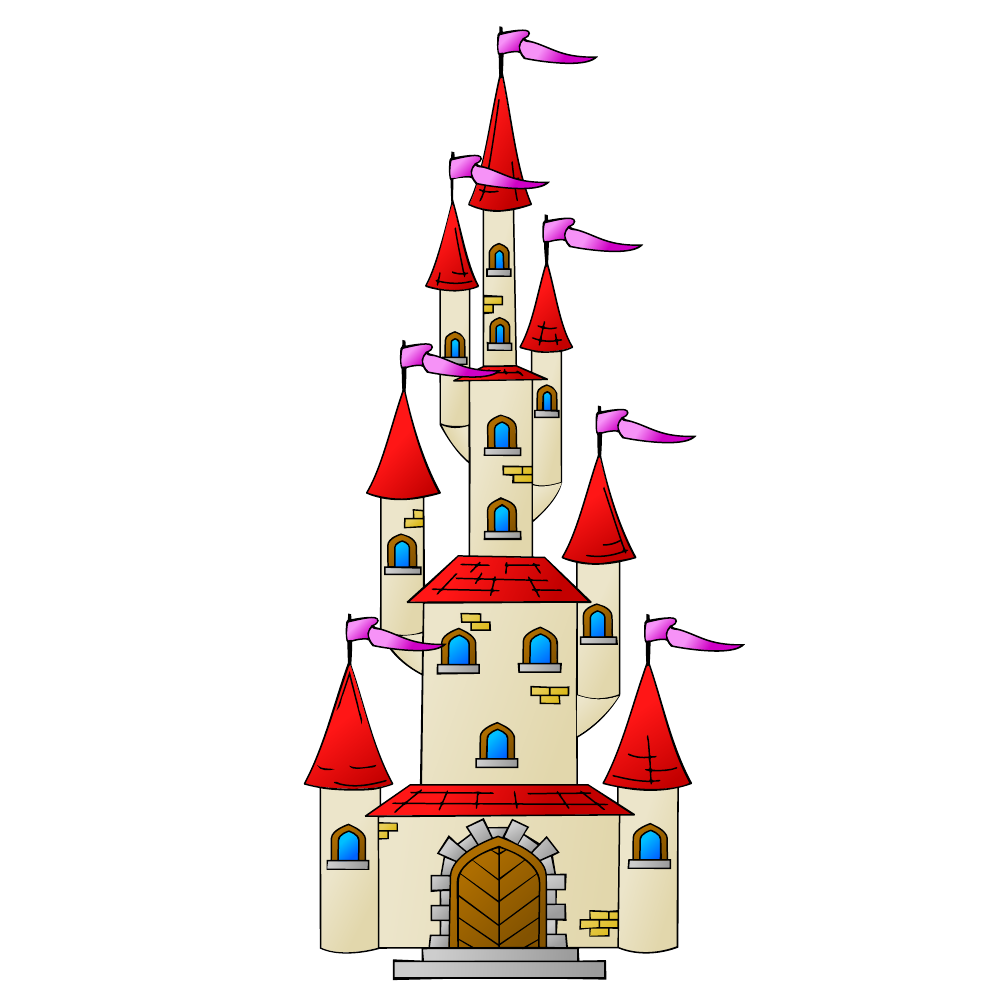}}; 
	\node (Music) at (-4,5) {\Dwarf{0.04}{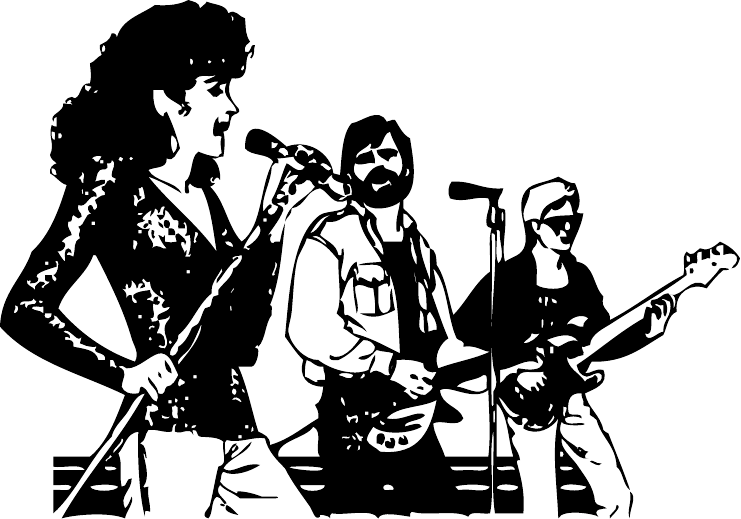}}; 
	\node (Truth) at (4,-5) {\Dwarf{0.03}{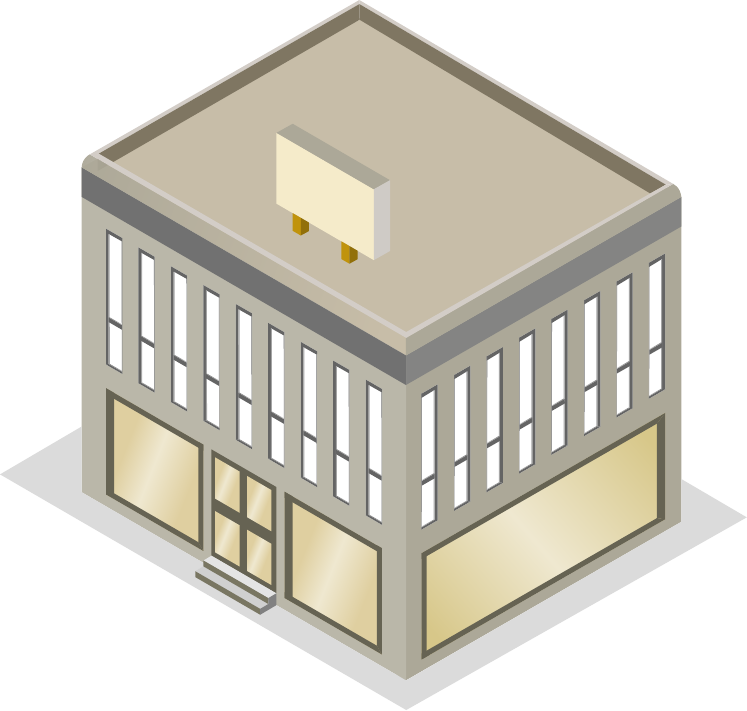}}; 
	\node (Math) at (-2.8,-4) {\Dwarf{0.03}{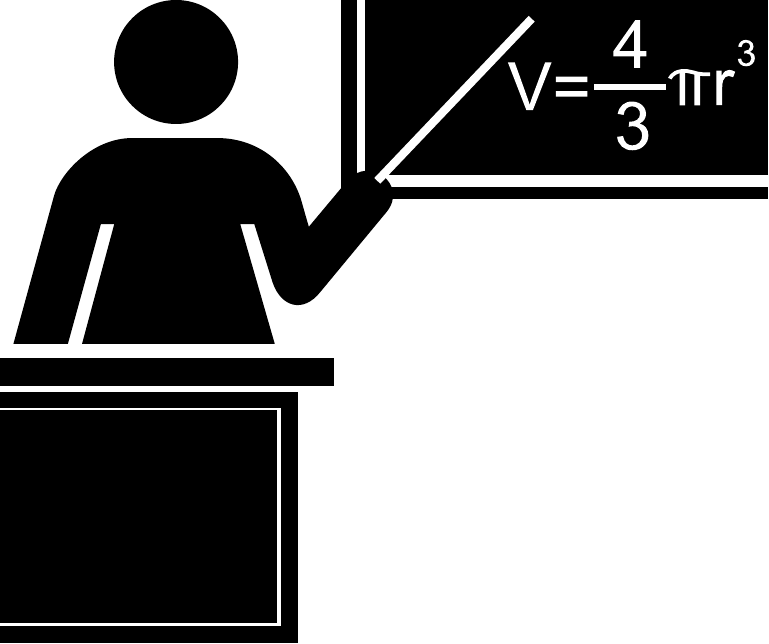}}; 
	\foreach \i  in {-6,...,6} {
		\fill [red] ($(\i,\i)$) circle (2pt);
	}
	\fill [green] ($(5,4)$) circle (2pt);
	\fill [blue] ($(2,3)$) circle (2pt);
	\end{tikzpicture}
	\caption{Dwarf Kingdom and Coordinate Geometry.} \label{fig:Dwarf_kingdom}
\end{figure}

It can be hoped, given the preceding discussion, that storytelling can also be an effective pedagogical tool in the teaching of mathematics, especially to young students. 
One advantage of an effective story is that it engages young curious students. Engagement is the first step towards learning. Indeed, it was a story of a `Dwarf Kingdom' that helped six-year-old Sophie in learning the basics of coordinate geometry. There was a land of dwarfs, where these little creatures moved by hopping. In one jump they could only hop one unit to their left, or right, or up, or down (see Fig.~$ \ref{fig:Dwarf_kingdom} $). In the center of the kingdom was the magnificent Palace of the King. A little small but equally splendid was the palace of the dwarf princess. There were other important buildings of the kingdom, such as the Ministry of Music, the Ministry of Magic, the Ministry of Truth, the Ministry of Beauty and the Ministry of Math. One fine day the king of dwarfs received a secret message from his trusted spy that the neighboring Kingdom of Elves is planning an attack on the kingdom of dwarfs.  Elves were very skilled magicians. Therefore, in order to prepare for their defense, the king of dwarfs ordered his messenger Neo to go from the King's palace to the Ministry of Magic, and urgently summon the royal magician. Now Sophie was asked to help Neo in reaching to the Ministry of Magic as quickly as possible. With Sophie's advice, Neo first hopped $ 3 $ times to his right and then $ 4  $ times up. Will he reach the correct destination? Of course, said Sophie. She had advised Neo correctly. Later she learned to write ($ 3 $ right , $ 4 $ up) to describe the motion of Neo in short, where it was agreed that `left or right' movements will be written first and thereafter 'up or down'. For example, ($ 3 $ left, $ 4 $ down) was also a valid move and if Neo followed this starting from the King`s palace, he will reach the Ministry of Math. It was agreed that hops to the left for $ n $ times  will be represented as $ -n $ hops to the right. Similarly, hops to $ -n $ times up represented $ n $ hops downwards. Therefore, ($ 3 $ left, $ 4 $ down) is the same as ($ -3 $ right, $ -4 $ up). This formulation also helped Sophie in reinforcing her understanding of negative numbers. 

Now with the impending danger of war, the royal magician created `a magical oil'. If this magical oil is sprinkled on a particular place, then any enemy soldier entering there will turn into a stone. Suppose the royal magician wanted to sprinkle his magical oil on the blue dot in Fig.~$ \ref{fig:Dwarf_kingdom}$ and Neo was again assigned this task. He had to start from the King's Palace with the magical oil and very carefully hop to his destination. Sophie once again helped Neo and wrote down the coded instruction for Neo as ($ 2 $ right, $ 3 $ up). Later on, more details were added to the story. In fact, it was turned into a small game. First, the condition was imposed that Neo is allowed to reach only those places where he can reach by hopping an equal number of times to right and up. Now an enemy position was given, say the green dot in Fig.~$ \ref{fig:Dwarf_kingdom}$. Sophie found out that Neo couldn't reach there under the constraints imposed on the movement of Neo. The game continued for a while. Sophie began to notice some patterns. Then Sophie was advised to list down all the places Neo could reach, if he was instructed to hop an equal number of times to the right and the up. She figured it out with some help and then marked all such places with red dots (see Fig.~$ \ref{fig:Dwarf_kingdom}$). She immediately noticed that all the red dots appeared to lie on a straight line. Of curse, Dwarf moved only on an integer lattice still Sophie was learning many interesting new things, such as the connection between the constraint on hopping `right $ = $ up' with the pattern of red dots appearing to lie on a straight line ($ y=x$). She even learned to plot all the points corresponding to the constraint on hopping `right $ > $ up'. It was a story, play, and a hands-on painting activity for her. She not only enjoyed this activity, but using her rich imagination she went on to add more to the story of the dwarfs and in the process learned some of the fundamental concepts of coordinate geometry.

Of course, the art of storytelling and its theory and techniques in the literary context is not  our concern here. We are interested in exploring how can storytelling be effectively used as a tool in teaching mathematics in a classroom. We can gain some insights on this from the example of the story of the Dwarf Kingdom. Young children are especially fond of fairy tales, and stories with monsters or animals as characters in them. It is a good idea to include such characters in our mathematical stories. A story should have some conflict or tension which can evoke emotions in the listener.  Moreover, if possible adding some dramatic words, sounds or actions while narrating the story is helpful. For example, in the above story of dwarfs, one can enact the hopping of dwarfs with a slouched gesture and perhaps a hilarious grunting sound after each hop. Such gestures and sounds are especially appealing to small kids and they tend to enjoy such stories. Moreover, a story in the educational context, and especially if a mathematics classroom, is even more effective if it leads to some hands-on activity. 

Mathematical Anecdotes, jokes on mathematicians or other forms of humor could also be effective in engaging students, especially older students. For example, a story in  \cite{barrow2008100}, describes how calculus literally saved the life of the Nobel prize winner physicist Igor Tamm, when he was caught by bandits during the Russian revolution period. 
\begin{tcolorbox}[breakable, enhanced, colback = mycolor]
	In the Russian revolutionary period, Tamm was a young professor teaching physics at the University of Odessa in the Ukraine. Food was in short supply in the city, and so he made a trip to a nearby village, which was under the apparent control of the communists, in an attempt to trade some silver spoons for something more edible, like chickens. Suddenly, the village was captured by an anti-communist bandit leader and his militia, armed with rifles and explosives. The bandits were suspicious of Tamm, who was dressed in city clothes, and took him to their leader, who demanded to know who he was and what he did. Tamm tried to explain that he was merely a university professor looking for food.
	`What kind of professor?' the bandit leader asked.
	`I teach mathematics,' Tamm replied.
	`Mathematics?' said the bandit. `All right! Then give me an estimate of the error one makes by cutting off Maclaurin's series at the $ n^{\text{th}} $ term. Do this and you will go free. Fail, and you will be shot!'
	Tamm was not a little astonished. At gunpoint, somewhat nervously, he managed to work out the answer to the problem - a tricky piece of mathematics that students are taught in their first course of calculus in a university degree course of mathematics. He showed it to the bandit leader, who perused it and declared `Correct! Go home!'
\end{tcolorbox}

Considered in the larger social context and as a collective human endeavor, mathematics needs role models and icons to motivate and inspire its practitioners.  The history of mathematics is a rich source of inspirational stories. Many present-day mathematicians have mentioned being inspired by past legends. Children should be introduced to biographies of geniuses such as Gauss, Euler, Riemann, Galois, Abel, Ramanujan, Erd\H{o}s etc. 

Another aspect of the history of mathematics is the fascinating story of how mathematical ideas have developed. A mathematical idea, which seems obvious to us, might have taken several hundreds of years to develop in its present simplified form.  Also, there are many important theorems and conjectures in mathematics, worth introducing to young students. For example, the statement of Fermat's last theorem is elementary and even school students can easily understand it.

\section{Meeting coordinate geometry early.} \label{Sect:EarlyCoordinate}
\begin{figure}
	\begin{subfigure}{.5\textwidth}
		\centering
		\includegraphics[scale=0.4]{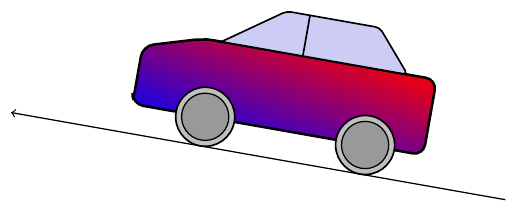}
		\caption{}\label{fig:car}
	\end{subfigure} 
	\begin{subfigure}{.5\textwidth}
		\centering
		\includegraphics[scale=0.4]{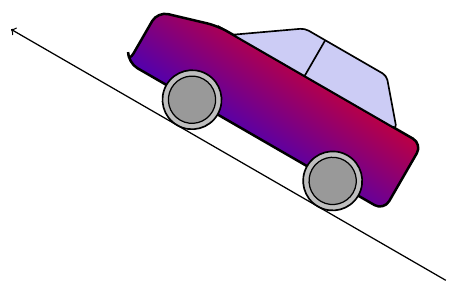}
		\caption{} \label{fig:carone}
	\end{subfigure}
	\caption{The concept of slope: a journey from the intuitive understanding of `slope' to the symbolic mathematical definition.}
\end{figure}

We have already described how coordinate geometry was introduced to Sophie using a story about the Dwarf Kingdom. Another method used was a game played on a big grid on the floor. Some small toys were lying on the floor at some of the grid points. A player gets the toy by correctly finding the coordinates of the point where a toy is located. 
One important concept in coordinate geometry is that of the slope of a straight line. In fact, even a six-year-old can easily figure out that the road in Fig.~$ \ref{fig:carone} $ is steeper than that of Fig.~$ \ref{fig:car} $. It is slightly harder to connect this intuition to the traditional definition of the slope as the ratio of `rise' over `run'. Still, after some initial difficulty, Arya understood the mathematical definition of the slope of a straight line. Arya also reviewed the basic concepts of equivalent fractions while learning about the slope. The concept of similar triangles was also introduced to him. Now by connecting all these different aspects related to the concept of slope,  his understanding of fractions was reinforced. He was also getting well-prepared for learning calculus later, where the understanding of the concept of the slope will be useful.

In Sect.~$ \ref{Sect:Physics_Eng} $ and Sect.~$\ref{Sect:Computer_games} $ we will again discuss some other strategies for teaching coordinate geometry to young gifted students.

\section{Connections: combining it all.} \label{Sect:Combine}
One source of aesthetic beauty in mathematics is the unexpected connections between apparently disconnected sub-fields. The joy of true understanding overwhelms a scientist or a  mathematician when she uncovers the hidden truth behind a phenomenon or a mathematical problem.  It is true that $ x=y $ is an algebraic relation, but at the same time, it also represents a geometric object, namely a straight line passing through the origin and having the slope $ 1 $. We have already emphasized the need to introduce algebra and coordinate geometry as early as possible to young gifted children. Also in Sect.~$ \ref{Sect:Distributive} $, we have noted that there exists an easy pathway for proving the Pythagorean theorem using the area model and the distributive property. Once this theorem is understood by a gifted child, several avenues open up for learning many interesting mathematical results. We note a few possible directions:
\begin{enumerate}
	\item Assuming that coordinate geometry is already introduced, students can now learn the distance formula, i.e., the shortest distance between the points $ (x_1,y_1) $ and $ (x_2,y_2) $ is given by $ \displaystyle{\sqrt{(x_2-x_1)^2 + (y_2-y_1)^2 }} $.
	\item An equation of a circle, say with center at $ (a,b) $ and radius $ r$, can be introduced as this equation $ \displaystyle{(x-a)^2 + (y-b)^2 =r^2}  $ depends on the distance formula or essentially on the Pythagorean theorem.
	\item The trigonometric identities such as $ \displaystyle{\sin^2 \theta + \cos^2 \theta =1 } $, $ \displaystyle{ 1  +  \tan^2 \theta = \sec^2 \theta} $ etc. are just alternate ways of expressing the Pythagorean theorem. Therefore, these concepts can easily be introduced quickly after the Pythagorean theorem is taught.
	\item Pythagorean theorem is used in finding the magnitude of a vector and also a complex number. 
\end{enumerate}
The above discussion makes it clear that a gifted student can benefit from an integrated learning approach that combines concepts from algebra, coordinate geometry, trigonometry, vectors and complex number together.  There is no reason that a gifted student is forced to wait for several years before getting a taste of these subjects. Unfortunately, due to excessive focus on arithmetic, this happens quite often. Of course, we are not suggesting that in one go, a gifted student will completely learn and master all the subjects mentioned earlier. What we are suggesting is that if opportunities open up for introducing advanced concepts while teaching, then without hesitation such advanced concepts should be introduced to gifted students. In fact, the learning would take place non-linearly, as it mostly happens for research mathematicians. 

\begin{figure}
	\mbs
	\centering
	\includegraphics[scale=0.6]{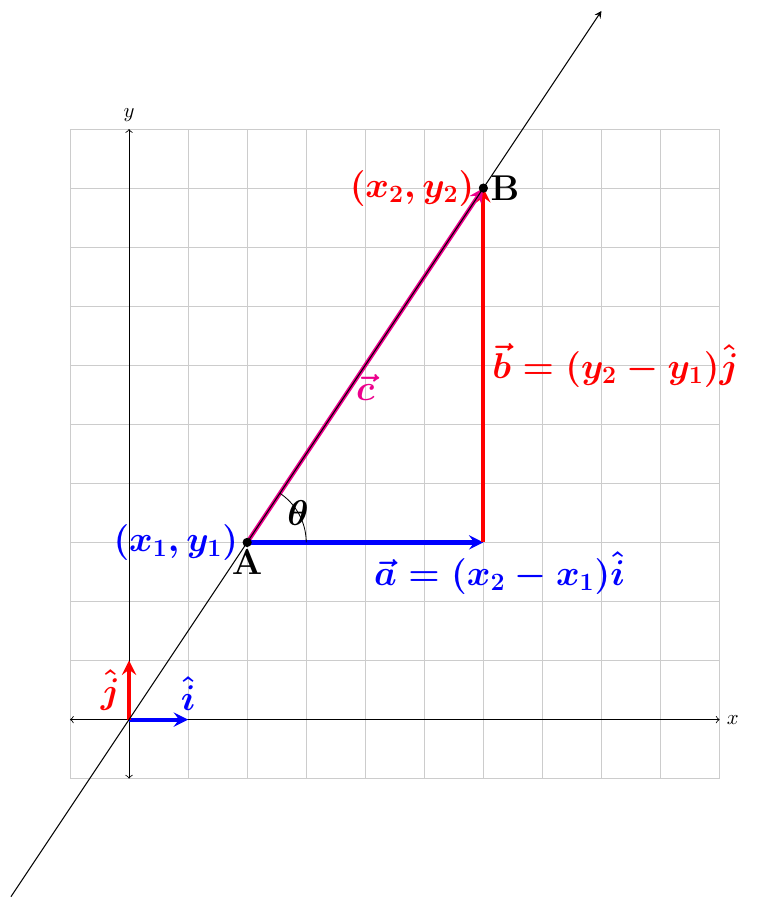}
	\mbe
	\caption{One diagram, connected with multiple concepts.}\label{fig:vectorcombined}
\end{figure}

We have remarked earlier that mathematics is highly interconnected, and these interconnections should be utilized in teaching gifted students. For example, from a single diagram (see Fig.~$ \ref{fig:vectorcombined} $), Pragya learned a number of mathematical concepts. 

\begin{enumerate}
	\item The triangular law of vector addition: $ \vec{c} = \vec{a} + \vec{b} $.
	\item Decomposition of a vector into its components:
	$$\vec{c} = (x_2-x_1) \hat{i} + (y_2-y_1)\hat{j} .$$
	\item Pythagorean theorem (or the distance formula):
	$$ |\vec{c}| = \sqrt{(x_2-x_1)^2 + (y_2-y_1)^2} .$$
	\item Describing a Vector with the given end points: The vector from point A with coordinates $ (x_1,y_1) $ to point B with coordinates 
	$ (x_2,y_2) $ is given by $$(x_2-x_1) \hat{i} + (y_2-y_1)\hat{j}  .$$
	\item The polar decomposition of a vector: $ \vec{a} = |\vec{c}| \cos \theta \, \hat{i}$ and   $\vec{b} =  |\vec{c}| \sin \theta \, \hat{j} .$
	\item The slope of a given line: The slope of the line joining points A and B is  $$ {m = \frac{\text{Rise}}{\text{Run}} = \tan \theta =  \frac{y_2-y_1}{x_2-x_1}} .$$
	\item An equation of the straight line: An equation of the straight line joining points A and B is $$ \displaystyle  {y - y_1 =  \left( \frac{y_2-y_1}{x_2-x_1} \right) (x-x_1) }.$$
\end{enumerate}

\section{Fractals.} \label{Sect:Fractals}
Benoit Mandelbrot coined the term Fractals in 1975 and introduced a class of mathematical objects,  that are beautiful and also appear in nature. Indeed, fractals are connected to a wide variety of natural objects such as the shapes of clouds, coastlines,  ripples in oceans and even flowers and trees. Of course, fractals produce stunning visual images and are used in producing realistic special effects in movies and video games. Fractals are defined by iteration of simple mathematical equations and it is indeed amazing that such simple mathematical relations produce complex fractal images. For instance, consider the iteration of the equation $ f(z) = z^2 + c $, starting from $ z=0 $ one can obtain the sequence $ z_0,\,z_1, \, z_2,\ldots $  with $ z_0 = f(0) =c $, $ z_1= f(z_0) = f(c) = c^2 +c $, $ z_2 = f(z_1) = (c^2 +c)^2 +c  $ and so on. The Mandelbrot set consists of all of those complex numbers $ c $ for which the elements of the sequence  $ z_0,\, z_1, \, z_2,\ldots $ remain bounded in absolute value. Fig.~$ \ref{fig:Mandelbrot}$ shows the picture of the Mandelbrot set. 

\begin{figure}
	\mbs
	\centering
	\includegraphics[scale=0.8]{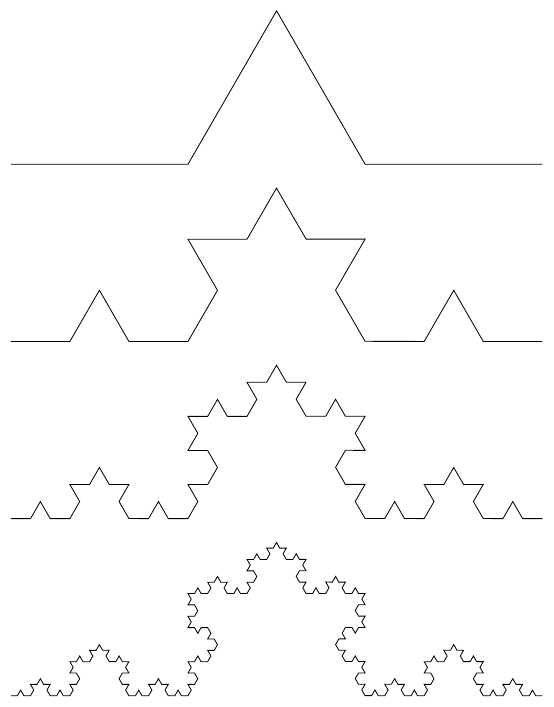}
	\mbe
	\caption{Koch snowflake.}	\label{fig:Koch_snowflake}
\end{figure}

When Pragya first encountered fractals, she was fascinated by their beauty. First, she learned how to draw Koch snowflakes 
(see Fig.~$ \ref{fig:Koch_snowflake} $) 
by hand. Then, as she was already familiar with the concept of recursion and the knew how to program in python using its library `turtle', it was a good opportunity for her to use her programming skills to draw some beautiful pictures of fractals (see Sect.~$ \ref{Sect:Scratch} $).

\begin{figure}
	\mbs
	\centering
	\begin{subfigure}{.3\textwidth}
		\centering
		\includegraphics[scale=0.5]{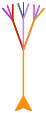}
		\caption{} \label{tree_50}
	\end{subfigure}
	\begin{subfigure}{.3\textwidth}
		\centering  	
		\includegraphics [scale=0.5] {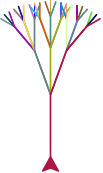}
		\caption{} \label{tree_100}
	\end{subfigure}
	\begin{subfigure}{.3\textwidth}
		\centering  	
		\includegraphics [scale=0.5]{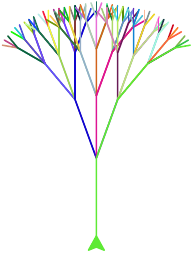}	
		\caption{} \label{tree_40}
	\end{subfigure}
	\begin{subfigure}{.5\textwidth}
		\centering  	
		\includegraphics [scale=0.2]{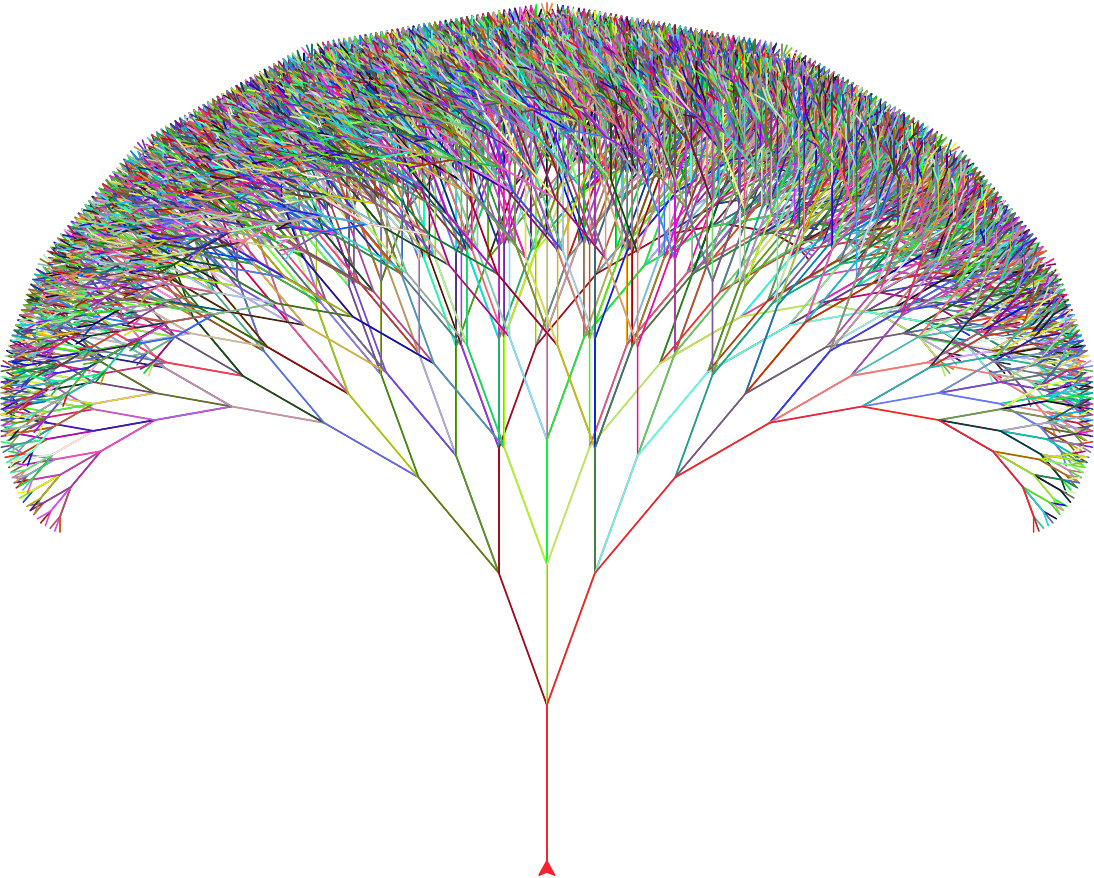}	
		\caption{} \label{tree_50}
	\end{subfigure}
	\mbe
	\caption{Computer-generated trees using recursion. The trees are drawn using the Turtle graphics library of Python with an increasing depth of recursion from Subfig.~$ (a) $ to Subfig.~$ (d) $.}  \label{fig:tree_recursive}
\end{figure}

\pgfdeclarelindenmayersystem{Koch curve}{
	\rule{F -> F-F++F-F}}
\pgfdeclarelindenmayersystem{Sierpinski triangle}{
	\rule{F -> G-F-G}
	\rule{G -> F+G+F}}
\pgfdeclarelindenmayersystem{Fractal plant}{
	\rule{X -> F-[[X]+X]+F[+FX]-X}
	\rule{F -> FF}}
\pgfdeclarelindenmayersystem{Hilbert curve}{
	\rule{L -> +RF-LFL-FR+}
	\rule{R -> -LF+RFR+FL-}}

\begin{figure}
	\mbs
	\begin{subfigure}{.45\textwidth}
		\centering
		\includegraphics[scale=0.6]{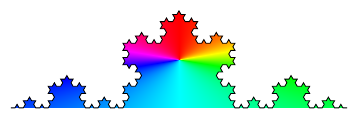}
		\caption{{\footnotesize \centering{Koch curve. (L-system: rule F $ \to $ F-F++F-F, angle=-60, axiom=F, order=4.)}}} \label{Koch}
	\end{subfigure}
	\begin{subfigure}{.45\textwidth}
		\centering
		\includegraphics[scale=0.3]{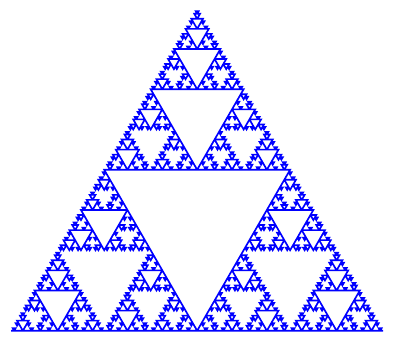}
		\caption{{\footnotesize Sierpinski triangle. (L-system: rule F $ \to $ G-F-G and G $ \to $ F+G+F, angle=-60, axiom=F, order=8.)}} \label{fig:Sierpinski}
	\end{subfigure}
	\begin{subfigure}{.45\textwidth}
		\centering
		\includegraphics[scale=0.45]{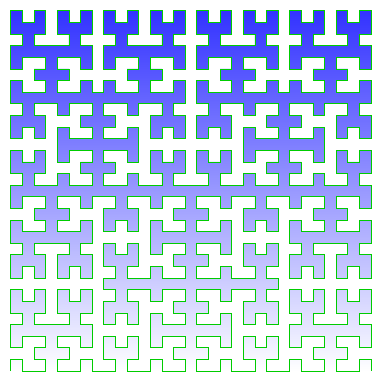}
		\caption{{\footnotesize Hilbert curve.(L-system: rule L $ \to $ +RF-LFL-FR+ and R $ \to $ -LF+RFR+FL-, angle=20, axiom=F, order=6.)}} \label{Hilbert curve}
	\end{subfigure} \hspace{1 cm}
	\begin{subfigure}{.45\textwidth}
		\centering
		\includegraphics[scale=0.45]{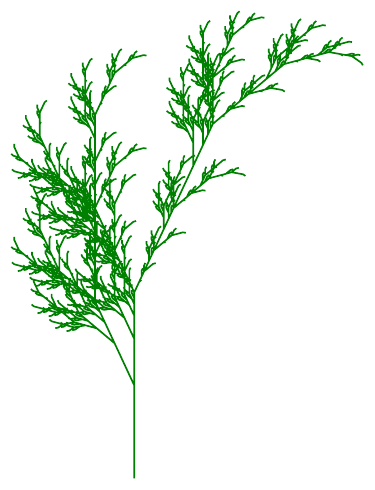}
		\caption{{\footnotesize Fractal plant. (L-system: rule X $ \to $ F-[[X]+X]+F[+FX]-X and F $ \to $ FF, angle=20, axiom=F, order=6.)}} \label{Fractal plant}
	\end{subfigure}
	\mbe
	\caption{Fractals using L-system.} \label{fig:L-system-all}
\end{figure}

\section{Introducing beautiful examples and results from advanced topics.} 
There is so much beautiful math that can be introduced early to gifted students. 
Of course, it will not be possible for gifted children to understand the proof of many such theorems, although sometimes intuitive explanations can be provided to them. 
Still, even knowing some of these beautiful results will be intriguing for them and it will kindle their interest in mathematics. 
For example, gifted students can be told about results from number theory and topology, such as Fermat's last theorem, Euclid's proof of the infinitude of primes, Twin prime conjecture, Collatz Conjecture, Jacobi's theta function and its connection with the period of a real pendulum, Euler's formula $ V-E+F=2 $ connecting vertices $ V $, edges $ E $, and faces $ F $ of complex polyhedrons, Brouwer fixed-point theorem and its illustration with the example of a coffee cup in which the coffee is stirred using a spoon, Borsuk-Ulam theorem illustrated as implying that there always exist a pair of antipodal points on the Earth's surface with the same temperatures and also the same barometric pressures.
\section{Experimenting, exploring and discovering.} 
At least some of the joy, that a research mathematician experiences while exploring the beautiful world of mathematics, should be accessible to gifted students.
One way to achieve this is to give students open-ended problems. We give a few examples in the following.
\begin{enumerate}[label=Ex.~\arabic*.]
	\item Find all the different patterns that you can find in Pascal's triangle. Can you guess what would be the sum of all the numbers in the $ n^{\textit{th}} $ row in the Pascal's triangle? (Here it is convenient to start with the topmost row being the $ 0^{\textit{th}} $ row).
	\item Find all the patterns or relations you can find in the sequence of Fibonacci numbers. A couple of examples are: $ \sum_{k=1}^{n} F_k = F_{n+2} -1 $ and $ \sum_{k=1}^{n} F_k^2 = F_n F_{n+1} $
	\item Determine the decimal expansion of $ 1/89 $ up to a few terms. Do you notice any patterns? Can you find the number in base $ 5 $ which has the same pattern?
	\item Can you find the Fibonacci Numbers in Pascal's Triangle? 
	\item Find out the heights of all the students in your class and then create a height-frequency plot. What is the shape of this graph? 
\end{enumerate}

\subsection{Interactive algebra and geometry software programs.} \label{Sect:Geogebra}

There are many interesting results in elementary mathematics, especially in euclidean geometry, that students can discover by doing experiments. Interactive Computer software tools such as \textbf{Geogebra} (\url{https://www.geogebra.org/}) can also be effectively used in these experiments. For example, in a guided inquiry-based learning exercise Medha was asked to draw two parallel lines and also a third transversal intersecting the two parallel lines. Then she was asked to measure various angles and after repeating this experiment several times she understood the relations between a pair of corresponding angles, alternate angles and consecutive interior angles. Similarly, she did experiments using Geogebra to measure all the internal angles of a triangle and found the sum to always be equal to $ 180 $ degrees. She was also guided to discover the Pythagorean theorem, the ratio of the circumference and the diameter of a circle and several other results in geometry by actual experimental measurements. This was a sort of active learning exercise for her. It helped in kindling her interest in geometry, and she was motivated to learn and understand the axiomatic proof-based geometry to understand the puzzling results that she had discovered through her experiments.

\begin{figure}
	\centering
	\includegraphics[scale=0.5]{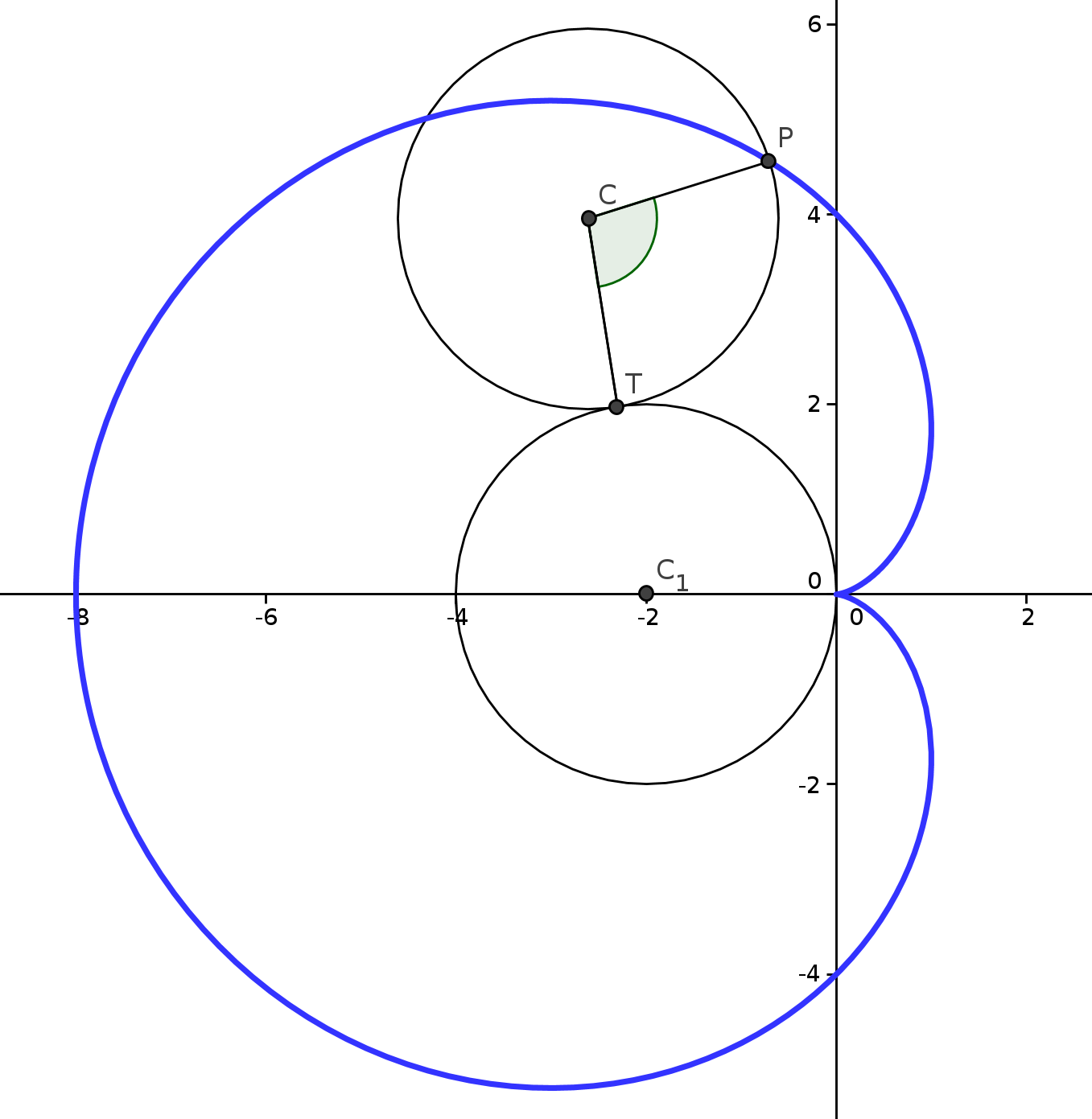}
	\caption{A cardioid. A cardioid is obtained by simulating the rotation of a circle, over the boundary of another circle of the same radius, using Geogebra. The locus of point P is a cardioid.} \label{fig:Cardiod}
\end{figure}

Pragya was given the following problem. Consider a special point on a wheel of radius $ r $ which is rotating without slipping on a straight path. A point P is marked on the wheel. Find the trajectory of point P.  Of course, the answer is that S traces out a cycloid as the wheel rotates and moves forward. She could really see this by doing a practical experiment with a wheel and later she simulated this using Geogebra. If instead, if we consider the trajectory of a point on the perimeter of a wheel, that is rolling around another fixed wheel of the same radius, then the resulting curve is a cardioid. Pragya simulated the rolling of one circle around the other and obtained the curve cardioid using Geogebra (see Fig.~$ \ref{fig:Cardiod} $). She was fascinated by observing this familiar curve appearing once again in a different context (see Sect.~$ \ref{Sect:Modular}$ for an earlier discussion).

\section{Puzzles, aha moments and creative problem-solving.} \label{Sect:Puzzle}
Recently in a discussion on mathematics in the lobby of the math department at the University of Manitoba a professor remarked - `` There are several types of mathematicians, with different personalities and tastes in life, some may be nice and humble, some arrogant, however, invariably every one of them must have experienced the `aha moment' of clarity of a concept or the joy of solving a hard problem.``
Indeed, one of the goals of early mathematics education should be that students experience this `aha moment'. Puzzles are important in the sense that they offer opportunities for young gifted students to experience the joy of problem-solving. However, it is important to realize that puzzles can also be extremely frustrating to students. Therefore, it is important to select good puzzles. Moreover, in case a student is totally stuck on a puzzle, appropriate hints and guidance should be provided. We also remark that a good puzzle for young children should be simple to state, easy to understand and fun to try. A puzzle that appears simple and solvable is less likely to intimidate a young child. Puzzles that involve actual physical activity such as tracing a curve or rearranging dominoes on a board, etc. hold special appeal to young children. 

In the following, we provide some illustrative examples.

\begin{figure}
	\mbs
	\begin{subfigure}{.45\textwidth}
		\centering
		\includegraphics[scale=0.5]{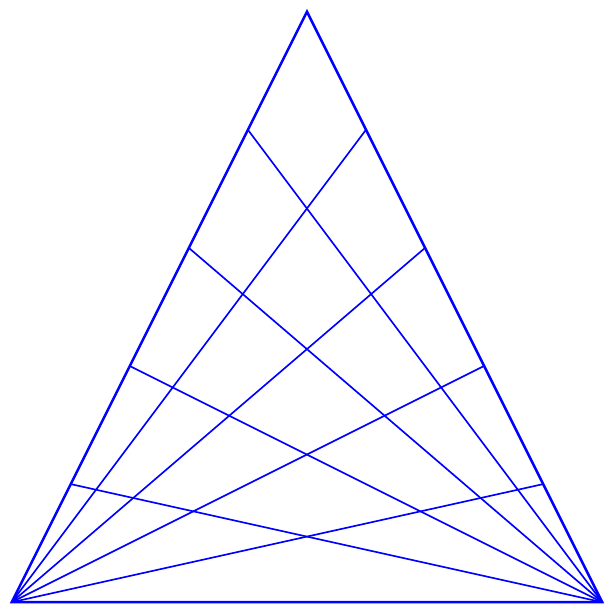}
		\caption{The equal sides of the triangle are divided into $ 5 $ segments.} \label{fig:triangle}
	\end{subfigure} \hspace{1 cm}
	\begin{subfigure}{.45\textwidth}
		\centering
		\includegraphics[scale=0.5]{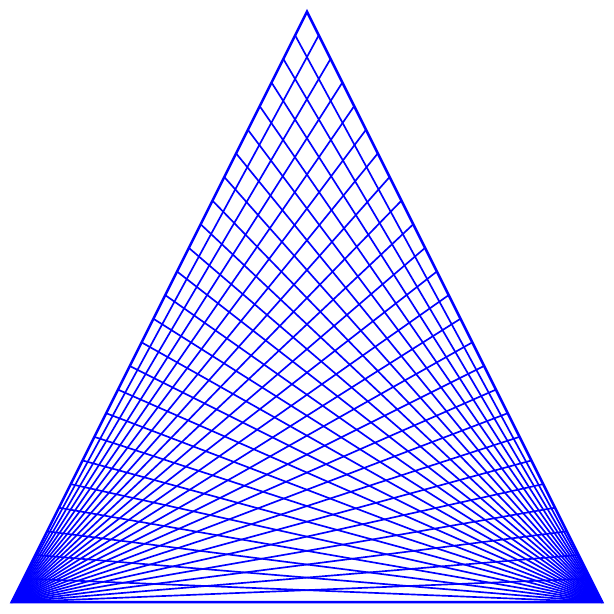}
		\caption{The equal sides of the triangle are divided into $ 25 $ segments. } \label{fig:triangle_25}
	\end{subfigure}
	\mbe
	\caption{How many triangles are there? }
\end{figure}

\begin{figure}
	\begin{subfigure}{.5\textwidth}
		\centering
		\includegraphics[scale=0.2]{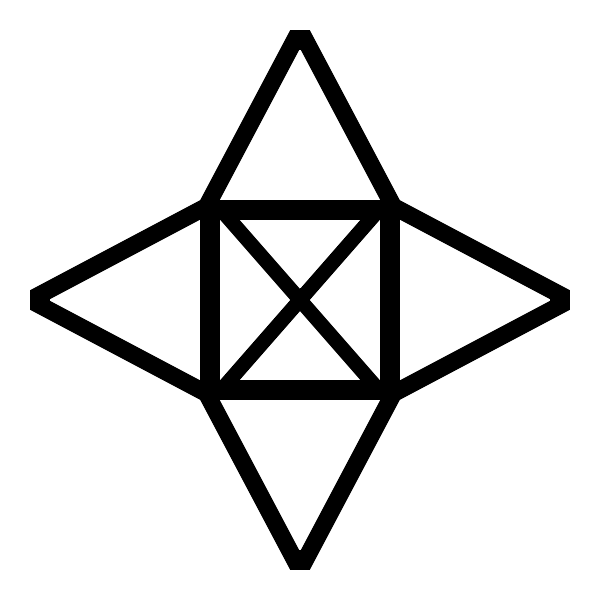}
		\caption{{\footnotesize Can you draw this shape without lifting your pen and without retracing any edge?}} \label{fig:graphtrace}
	\end{subfigure}
	\begin{subfigure}{.5\textwidth}
		\centering
		\includegraphics[scale=0.4]{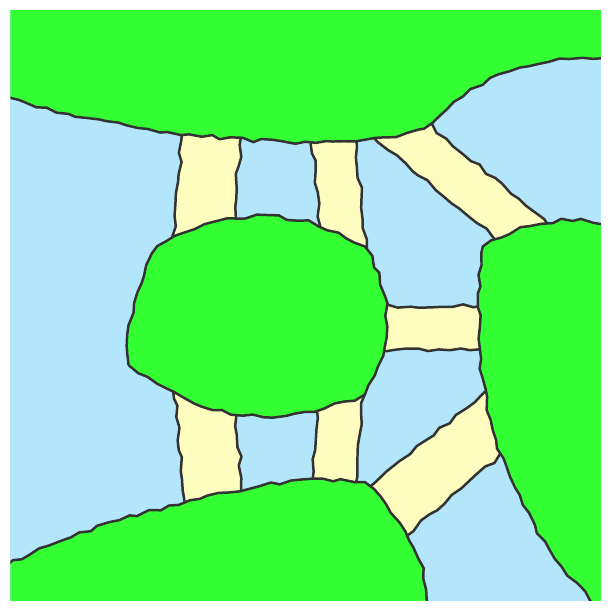}
		\caption{{\footnotesize The Seven Bridges of K\"{o}nigsberg problem. }} \label{fig:Konigs}
	\end{subfigure}
	\caption{Puzzles related to Eulerian paths.}
\end{figure}

\begin{enumerate}[label=Ex.~\arabic*.]
	\item \textbf{Counting the number of squares on a chess board.} When Arya was asked how many squares are there on a regular $8 \times 8  $ chess board, he immediately replied $ 64 $. However, when he was told that there are more squares, for example, there is a big square of the size $ 8 \times 8 $, he got stuck on the problem. Arya was provided with some hints. One important skill in mathematical problem solving is to start with a smaller or simpler baby version of the problem. Arya could solve the problem of counting the number of squares on a $ 3 \times 3 $ chess board. Then he went on to solve the problem for  $ 4 \times 4 $ chess board. Equipped with the experience gained from the smaller examples, he could finally solve the original problem for the $ 8 \times 8 $ chess board.
	
	\item \textbf{Drawing shapes in one stroke without retracing.} 
	Pragya was given the task of drawing the diagram in Fig.~$ \ref{fig:graphtrace} $, without lifting the pen, and without tracing the same line more than once. She believed that she could do it and therefore she tried for a long time. When she was getting frustrated, she was given another variant of the same problem, which she could solve. She discovered that some of the shapes are easy to draw under the given restrictions, while some others are hard or perhaps impossible. This puzzle provided a good opportunity for her to learn about Eulerian paths and the general result given by Euler that tracing a graph without lifting the pen and without repeating the same edge twice is possible only if  either all, or all but two vertices are of even degrees. Here the degree of a vertex is the number of vertices with which it is connected.
	\label{Enum:drawing}
	
	\item \textbf{The Seven Bridges of K\"{o}nigsberg.} Pragya was introduced to the famous problem related to the seven bridges of K\"{o}nigsberg. A sketch of the map of the city is given in Fig.~$\ref{fig:Konigs}$. The problem was to devise a walk through the city that would cross each of those bridges once and only once. Since she had already attempted the previous problem, with a little guidance (that she can represent the problem by a graph by viewing the land masses as vertices and the bridges as edges of a graph), she discovered that it was impossible to devise a desired walk.
	\label{Enum:Konigsberg}
	
	\item \textbf{Triangle counting problem.}  How many triangles are there in Fig.~$ \ref{fig:triangle} $ and Fig.~$ \ref{fig:triangle_25} $. 
	
	\item \textbf{Covering a square with a missing corner with dominoes.}  How can a $\bf 4 \times 4  $ square, whose two opposite corners are removed, be tiled using seven dominoes? One domino is  of the size of $ 2 \times 1 $ (see Fig.~$\ref{fig:Dominos}$). What about an $\bf 8 \times 8 $ square with one pair of opposite corners missing?
	
	\item \textbf{Sudoku for kids.} For very young kids the smaller $ 4 \times 4 $ version of Sudoku can be given (see Fig.~$ \ref{fig:Sudoku} $). 
	\item \textbf{Ants on a scale.}       Fifty ants are moving on a one-meter scale, which is so thin that an ant can travel only to the left or to the right. Each ant is traveling at a constant speed of $ 1 $ meter per minute. When two ants coming from the opposite direction collide, they simply reverse their directions of travel. When any ant reaches an end of the scale, it falls off.
	What is the minimum time after which it becomes certain, irrespective of the initial configurations of the ants on the scale, that all the ants have fallen off the scale?
\end{enumerate}

We remark that some of the puzzles that we presented have no solutions, while some others have multiple solutions. This is another important thing that students should be made aware of, to give them a sense of what lies ahead if they happen to study advanced mathematics or do research in this subject. Otherwise, most of the students erroneously assume that there is always a correct answer to any given math problem, such as the ones they find in their math exams, and mathematics is all about somehow getting to this correct answer. Once they encounter such problems they realize that there is more to mathematics than just following some algorithm to get to the right answer.

The students should be encouraged to think out of the box and ask creative questions. Can $ \sqrt{-1} $ be a member of our system of numbers? What can go wrong? Can we have a system where the multiplication operation is not commutative? Can we have a number system where $ 2 + 2 =5 $ and $ 1 + 1 =11 $? See \cite{padmanabhan_shukla_2022} for an interesting discussion on this. Asking such questions enables students to see mathematics as a creative discipline rather than a dull boring subject involving mindless calculations.

\begin{figure}
	\mbs
	\begin{subfigure}{.5\textwidth}
		\centering
		\begin{shidoku}
			\setrow{1}{1,2,3,4}
			\setrow{2}{3,4,1,2}
			\setrow{3}{4,3,2,1}
			\setrow{4}{2,1,4,3}
		\end{shidoku}
		\caption{Sudoku for kids.} \label{fig:Sudoku}
	\end{subfigure}
	\begin{subfigure}{.5\textwidth}
		\centering
		\begin{tikzpicture}[scale = 0.8]
		\draw (1,0) -- (4,0);
		\draw (0,4) -- (3,4);
		\draw [dashed] (0,0) --(1,0);
		\draw [dashed] (0,0) --(0,1);
		\foreach \y in {1,2,3}
		{
			\draw (0,\y) -- (4,\y);
		};
		\draw (0,1) -- (0,4);
		\draw (4,0) -- (4,3);
		\draw [dashed] (3,4) --(4,4);
		\draw [dashed] (4,3) --(4,4);
		\foreach \x in {1,2,3}
		{
			\draw (\x,0) -- (\x,4);
		};
		\foreach \x in {1,2,3}
		{
			\draw (\x,-0.5) -- (\x,-1.5);
		};
		\draw (1,-0.5)--(3,-0.5);
		\draw (1,-1.5)--(3,-1.5);
		\node (A) at (2.1,-2) {Domino};
		\end{tikzpicture}
		\caption{A $ 4 \times 4 $ square with a pair of diagonally opposite squares removed. Can it be completely covered by $ 2 \times 1 $ dominos? What about the same problem with $ n \times n $ square?} \label{fig:Dominos}
	\end{subfigure}
	\mbe
	\caption{Puzzles.}
\end{figure}


A typical math textbook contains problems at the end of each of its chapters. Many such problems are routine problems,i.e., problems which can be solved by a direct application of the  theory discussed in the concerned chapter. However, in general, solving non-routine mathematical problems, such as some of the puzzles discussed earlier, is not so easy even for a gifted student, as such problems often demand out-of-box thinking and creativity. Still, some basic principles can help students develop their toolbox for solving such creative problems. 
In his excellent book \textit{How To Solve It} \cite{polya2004solve}, Polya has given four basic principles of problem-solving: Understand the problem, Devise a plan, Carry out the plan, Look back; and he has discussed these principles in detail. Of course, these principles are a good starting point. Eventually, in the course of their mathematical journey students should develop, from their individual experiences, personalized toolboxes of mathematical knowledge, principles, and problem solving tricks and techniques.  
We list below some of the commonly used techniques or principles, that are helpful in mathematical problem-solving. 

\begin{itemize}
	\item  Pigeon-hole principle.
	\item  Inclusion-exclusion principle.
	\item  Method of distinguished element.
	\item  Double counting.
	\item  Mathematical induction.
	\item  Divide and conquer.
	\item  Recurrence relation. 
	\item  Principle of invariants. 
	\item  Principle of symmetry.
\end{itemize}

Here, we remark that sometimes a gifted student, who is not used to struggling with routine mathematical problems, may find the experience of solving a hard mathematical problem very frustrating. It is possible that such a student may even lose interest in mathematics. Therefore, dealing with failure and managing emotions in such a student is of great importance. We have further discussed this in Sect.~$ \ref{Sect:Failure} $.

\section{Mathematics and programming} \label{Sect:programming_physics}
One of the key themes of this work is to strongly argue in favor of integrating computer programming with mathematics education of young gifted children. Learning to program a computer is an exceedingly important skill in the modern technology driven world. Indeed, computer simulation is pretty useful in mathematical modeling, not only in natural sciences and engineering, but also in economics, geography, psychology and many other disciplines in humanities. It may appear that learning to program is a skill that is too advanced for a child to learn. However, this fear is misplaced and it is certainly possible for a gifted child to learn the fundamentals of computer programming. Moreover, there are many resources available for a child to learn to program a computer. We will discuss some of them, such as Scratch and Turtle graphics library in the following. Once a child learns the basics of programming, it offers a great tool to experiment and play with various ideas in mathematics and thereby providing a very hands-on way of learning for a child. Further, in many situations, it allows a gifted child to have concrete realizations of  abstract mathematical concepts.

\subsection{{Scratch, Python, Turtle and L-systems.}} \label{Sect:Scratch}
Scratch (see \url{https://scratch.mit.edu/}), developed by MIT Media Lab, is a visual programming language that is especially useful in teaching a child to code. Scratch programs are created using blocks that can be dragged and dropped. Students can use Scratch to create and share stories, animations and games and thereby expressing their creativity and imagination. Scratch offers a perfect logical foundation for students and after learning Scratch they can comfortably move on to learn other traditional mainstream programming languages.  

\begin{figure}
	\begin{subfigure}{.45\textwidth}
		\centering
		\begin{python}
			from turtle import *
			
			""" Draw a square of size 200 """		
			forward(200)
			left(90)
			forward(200)
			left(90)
			forward(200)
			left(90)
			forward(200)
			left(90)
			done()
		\end{python}
		\caption{A Python program to draw a square of the size 200 using Turtle graphics library.} \label{fig:square}
	\end{subfigure} \hspace{1 cm}
	\begin{subfigure}{.45\textwidth}
		\centering
		\includegraphics[scale=0.5]{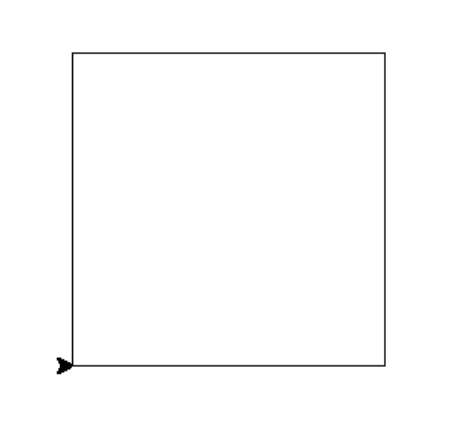}
		\caption{A square of size 200 was drawn using the Turtle graphics library  of Python.}
	\end{subfigure}
	\caption{A simple Turtle graphics program to draw a square.}
\end{figure}

Another good way to introduce programming to children is by using Turtle graphics, which like Scratch provides visual feedback. In fact, programming with quick visual feedback is extremely helpful for beginner programmers. Turtle graphics has its roots in the programming language Logo created by Wally Feurzeig, Seymour Papert, and Cynthia Solomon in 1967. Even a modern language like Python has support for Turtle graphics, which we will consider in the following. In Turtle graphics, a turtle can move relative to its current position by the following main commands.
\begin{enumerate}
	\item \textbf{Moving forward/backward.} Suppose the turtle is initially at $ (0,0) $ facing the direction of the positive $ x $-axis,  then turtle.forward(25) will move the turtle to $ (25,0) $ and also draw a line from $ (0,0) $ to $ (25,0) $ as the turtle moves. We note that turtle.forward(-25) or turtle.backward(25), both commands will result in the turtle moving to $  (-25,0)$. This can be used to reinforce the concept of negative numbers in a child.
	\item \textbf{Turning left/right.} Again assuming that the initial position of the turtle is at $ (0,0) $, pointing towards the positive $ x $-axis, we note that the command turtle.left(45) will turn the turtle by $ 45 $ degrees in place counterclockwise (i.e., in the left direction relative to its current direction of positive $ x $-axis).	
\end{enumerate} 
There are some other commands. In fact, the commands in Turtle graphics are quite simple and intuitive. The motion of the turtle on the screen makes it easy for children to learn and play with Turtle graphics. 
\begin{tcolorbox}
	\textbf{A classroom exercise.} A volunteer student, say Adam, is acting as a turtle. He is following the commands just like a turtle in Turtle graphics. Other students should be given some drawing challenges. For example, Clara has to draw a rectangle of given dimensions. She will issue a sequence of commands to Adam, such as forward(5), left(90) etc. to move him on a rectangle as required in the drawing challenge discussed in Sect.~$ \ref{Sect:Puzzle} $.
\end{tcolorbox}
Moreover, children can learn basic programming concepts such as that of a 'for-loop' using Turtle graphics. For example, in Fig.~$ \ref{fig:square} $ a square of size $ 200 $ is drawn using the Turtle graphics in Python by turning forward by 200 units and then turning left by 90 degrees in succession four times. Medha learned to write the same program by using a for loop (see Fig.~$\ref{fig:square_loop}  $).

\begin{figure}
	\centering
	\begin{python}
		from turtle import *
		
		""" Draw a square of size 200 using  a for-loop  """		
		for i in range(4):
		forward(200)
		left(90)
		done()
	\end{python}
	\caption{A Python program using Turtle graphics library to draw a square of the size 200 using a `for-loop'.} \label{fig:square_loop}
\end{figure}
Turtle graphics is also useful in learning geometry. For instance, after drawing a square the next challenge for Medha was to draw a regular hexagon. Of course, in order to draw a regular hexagon, she needed to find the angle of turn which led to her excursion in geometry. She learned how to compute an internal angle of a regular hexagon and also for any regular polygon. Finally, she could write a program to draw a hexagon. 

Pragya was helping Natasha in the activity that we discussed in Sect.~$ \ref{Sect:Counting} $. She created the patterns Fig.~$ \ref{fig:sum10} $, using Turtle graphics, where a point lying on any of the two oblique lines was joined to a point on the horizontal line if their sums equaled  $ 100 $. 

Turtle graphics can be used to draw objects with recursive structures such as fractals. For example, Pragya created a Python program using Turtle graphics  (see  Fig.~$ \ref{fig:turtle_trees_code} $) to draw the trees in Fig.~$ \ref{fig:tree_recursive} $. 
\begin{figure}
	\centering 
	\begin{python}
		from turtle import *
		import random
		
		screen.colormode(255)
		
		""" Recursively draws a tree using a turtle object t. 
		Len: is the length of the tree.
		theta: is the angle of turn (in degrees) of a child branch from its parent branch."""
		
		def tree(Len,theta,t):
			if Len < 5:
				return
			else:
				color(random.randint(0,255),random.randint(0,255),random.randint(0,255))
				forward(Len); left(theta); tree(Len-10,theta,t);
				right(theta); tree(Len-10,theta,t);
				right(theta); tree(Len-10,theta,t);
				left(theta);  backward(Len);
			return
		left(90); up(); bk(200); down();
		tree(100,20,t);
	\end{python}
	\caption{Recursively draws a tree using Turtle graphics. The output of this program is given in Fig.~$ \ref{fig:tree_recursive} $.}	\label{fig:turtle_trees_code}
\end{figure} 

\begin{figure}
	\centering
	\includegraphics[scale=0.8,frame]{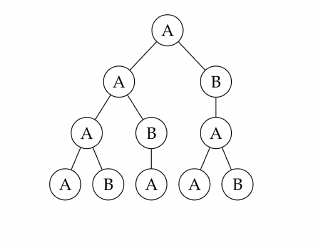}
	\caption{Evolution of the first few generations of an L-system which consists of the alphabets A and B, the axiom A and the rules A $ \to $ AB and B $\to  $ A. } \label{fig:L-system_AB}
\end{figure}

Turtle graphics can be combined with L-Systems to create  complex natural objects such as trees, flowers, snowflakes and many other organisms. An \textbf{L-system}, also known as the \textbf{Lindenmayer-system}, was named after the biologist and botanist Aristid Lindenmayer. An L-system consists of the following:

\begin{enumerate}
	\item \textbf{Alphabet.} A set of symbols that are used to form a word or string in the system. There are two types of alphabets. \begin{itemize}
		\item Constants: do not change during iterations. 
		\item Variables: change during iterations according to the production rules.
	\end{itemize}
	\item \textbf{Axiom.} An initial starting word from which the system grows.
	\item \textbf{Rules.} A set of rules that govern the growth of the system by stipulating how  the variables can be replaced with combinations of constants and other variables in each iteration. 
\end{enumerate}
For example, consider an L-system consisting of the alphabets A and B, the axiom A and the rules A $ \to $ AB and B $\to  $ A. This means in each iteration all the occurrences of A are replaced by AB and all the occurrences of B are replaced by A. 
We list below the growth of the system for the first few iterations (also see Fig.~$ \ref{fig:L-system_AB} $).
\begin{enumerate}
	\item \textbf{Axiom:}  A.
	\item \textbf{First iteration:} AB.
	\item \textbf{Second iteration:} ABA.
	\item \textbf{Third iteration:} ABAAB.
\end{enumerate}

In order to draw geometrical shapes, the variables in an L-system can be associated with the movements of a turtle. 
For example,
\begin{enumerate}
	\item \textbf{F} : means the turtle moves forward a certain distance (say d units) in the current direction.
	\item  \textbf{f}  : means the turtle moves forward a certain distance (say d units) in the current direction without drawing the line.
	\item  $ \mathrm{+} $ : means the turtle turns left by a certain angle.
	\item  $ \mathrm{-} $ : means the turtle turns right by a certain angle.
	\item \textbf{[ }: means save the current state of the turtle (i.e., its position and direction). It is like pushing the program state in a stack.
	\item \textbf{] }: means restore the last saved state of the turtle. It is like popping the program state from a stack. 
\end{enumerate}

Here we note that L-systems are supported by the TikZ and PGF Packages for Latex. All the sub-figures in  Fig.~$ \ref{fig:L-system-all} $ are drawn using the Lindenmayer System Drawing Library provided by TikZ and PGF.
For example, the following code 
\begin{tcolorbox}
	\begin{verbatim}
	\begin{tikzpicture}
	\draw [green!50!black, rotate=90]
	[l-system={rule set={{X -> F-[[X]+X]+F[+FX]-X},{F -> FF}}, 
	axiom=X, order=6, 
	step=2pt, angle=-25}]
	lindenmayer system;
	\end{tikzpicture}
	\end{verbatim}	
\end{tcolorbox}
generates the plant shown in Fig.~$ \ref{Fractal plant} $.

\subsection{{Computer graphics and image processing.}}
Computer graphics is full of applications of mathematics. Of course, most of these applications are too advanced for a child to learn. Still, there are plenty of attractive bits and pieces available in this field to make learning fun for a gifted child. 
We will discuss some of these ideas in the following. 
\begin{enumerate}[label=Ex.~\arabic*.]
	\item \textbf{Mathematical function describing a photograph.} With the easy availability of smartphones with a camera, it is easy to capture photographs of objects such as flowers, animals, buildings, clouds, water fountains, glass, windows, fruits etc. Then these objects can be imported into image editing software programs. In fact, in many of these objects students can see the occurrences of common mathematical curves, such as circles, ellipses, parabolas, cycloids etc. Sometimes a mathematical equation will not be apparent (but it can still be analyzed mathematically, see the remark at the end of this example). Still for advanced gifted students in high schools, these images can further be manipulated by importing them in a Python program. An image can be thought of a $ m \times n $ matrix with the entry at $ (i,j)^{\text{th}} $ position storing the integer defining the color at that position. Python supports many functions for working with an image. For example, consider the image of a cat as shown in Fig.~$\ref{fig:cat} $. The python program shown in Fig.~$\ref{fig:cat_code_rotate} $ extracts the contour of this image (which means the coordinates of the set of boundary points), see Fig.~$ \ref{fig:cat1} $, and then rotates the image counterclockwise by the angle $ \theta = \frac{\pi}{2} $, see Fig.~$ \ref{fig:cat2} $.

	\noindent \textbf{A remark.} We note that using FFT (Fast Fourier Transform) algorithm the set of boundary points can be transformed in the frequency domain. Then on carrying out the inverse FFT, one can get a function, which is a linear combination of complex exponential basis functions. Therefore, a complex function representing the shape of the contour can be obtained. Of course, FFT  is an extremely important algorithm in signal and image processing, but we do not expect it to be introduced to a school student. We just wanted to point out that for a given image (unless it is too pathological), one can obtain a function representing the shape of its contour.

	\begin{figure}
		\begin{subfigure}{.9\textwidth}
			\centering
			\begin{python}
				import math as mt
				import numpy as np
				import matplotlib.pyplot as plt
				from skimage import io
				from skimage import measure
				
				# Read the image of the cat
				cat = io.imread('cat.png')
				
				# Find the contour of the image cat
				contours = measure.find_contours(cat, 0.8)
				
				
				for n, contour in enumerate(contours):
				x = contour[:, 0]
				y = contour[:, 1]
				
				# Rotate the image by the angle theta 
				theta = mt.pi/2 
				x1 = x*mt.cos(theta)  - y*mt.sin(theta)
				y1 = x*mt.sin(theta)  + y*mt.cos(theta)
				
				plt.plot(x1,y1,linewidth=2)
				plt.show()
			\end{python}
			\caption{Extracting the contour of an image, see Fig.~$ \ref{fig:cat1} $, and then rotating the image counterclockwise by the angle $ \theta = \frac{\pi}{2} $, see Fig.~$ \ref{fig:cat2} $.} \label{fig:cat_code_rotate}
		\end{subfigure} \hspace{2 cm}

		\begin{subfigure}{.9\textwidth}
			\centering
			\begin{python}
				import numpy as np
				import matplotlib.pyplot as plt
				from skimage import io
				
				# Read the image of the cat
				camera = io.imread('cat.png')
				
				# Get the size of the image 
				(l_x, l_y) = camera.shape
				
				# Initialize the grid of size l_x times l_y
				X, Y = np.ogrid[:l_x, :l_y]
				
				# Compute which pixels should be masked. 
				outer_disk_mask = (X - l_x / 2)**2 + (Y - l_y / 2)**2 > (0.45*l_x)**2
				
				""" Apply the mask which will change the pixel color outside a disk 
				centered at (l_x/2,l_y/2) and of radius 0.45*l_x"""
				camera[outer_disk_mask] = 200
				
				plt.imshow(camera, cmap='gray', interpolation='nearest')
				plt.axis('off')
				plt.show()
			\end{python} 
			\caption{Creating a disk shaped mask using the image of the cat shown in Fig.~$ \ref{fig:cat}$.} \label{fig:cat_code_mask}
		\end{subfigure}
		\caption{Python programs for some simple image processing tasks. The programs use mathematical concepts like the rotation matrix for rotating the image, and the equation of the exterior of a disk to create a disk-shaped mask.  }
	\end{figure}

	\item \textbf{Transforming an image (translation and rotation).} Translation and rotation of images are required in many situations, for example, in a video game. Understanding of linear transformation is an important topic within mathematics itself. In fact, the rotation matrix helps in providing an intuitive understanding of various trigonometric identities and sum and difference formulas. It is advantageous to holistically teach interrelated topics such as the geometric meaning of multiplication by a complex number, Euler's formula for $ e^{i \theta} $ and the action of a rotation matrix on a vector. The multiplication action of a rotation matrix $ r(\theta) = \displaystyle{\mat{\cos \theta}{- \sin \theta}{\sin \theta}{\cos \theta} } $ on $ \vec{v} = \left[ \begin{array}{c}
	x
	\\ 
	y
	\end{array} \right]$ is to rotate the vector by the angle $ \theta $ in the counterclockwise direction, i.e., $ \vec{v_1} = r(\theta) \vec{v} $ is the vector obtained by rotating $ \vec{v} $ in the counterclockwise direction. For rotating an image by an angle $ \theta $, the rotation matrix needs to be applied for each of its points. 
	
	\begin{tcolorbox}
		In Fig.~$ \ref{fig:cat_code_rotate} $ the following lines of code 
		\begin{verbatim}
		x1 = x*mt.cos(theta)  - y*mt.sin(theta)
		y1 = x*mt.sin(theta)  + y*mt.cos(theta)
		\end{verbatim}
		are responsible for the rotation of $ 90 $ degrees of the image of the cat in Fig.~$\ref{fig:cat1} $ resulting in the image of the cat in  Fig.~$\ref{fig:cat2} $.  The above lines of code correspond to the following obvious mathematical calculation 
		\begin{align}
		\left[ \begin{array}{c}
		x_1
		\\ 
		y_1
		\end{array} \right] = \mat{\cos \theta}{- \sin \theta}{\sin \theta}{\cos \theta} 
		\left[ \begin{array}{c}
		x
		\\ 
		y
		\end{array} \right]
		\end{align}
		where the coordinate $ (x_1, y_1) $ is obtained by rotating the coordinate  of the point $ (x,y)$ counterclockwise by the angle $ \theta $. 
	\end{tcolorbox}

	\begin{figure}
		\mbs
		\begin{subfigure}{.5\textwidth}
			\centering
			\includegraphics[scale=0.45]{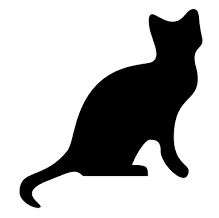}
			\caption{A cat.} \label{fig:cat}
		\end{subfigure}
		\begin{subfigure}{.5\textwidth}
			\centering
			\includegraphics[scale=0.4]{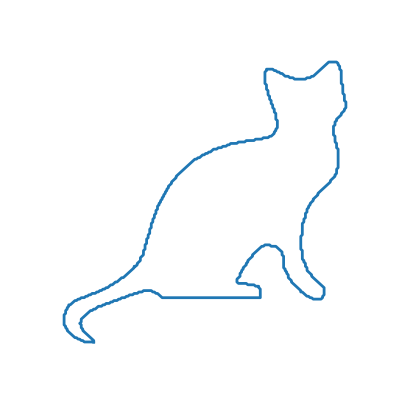}
			\caption{Extracting the boundary (contour) of the image of a cat using Scikit image processing library.} \label{fig:cat1}
		\end{subfigure}
		\begin{subfigure}{.5\textwidth}
			\centering
			\includegraphics[scale=0.4]{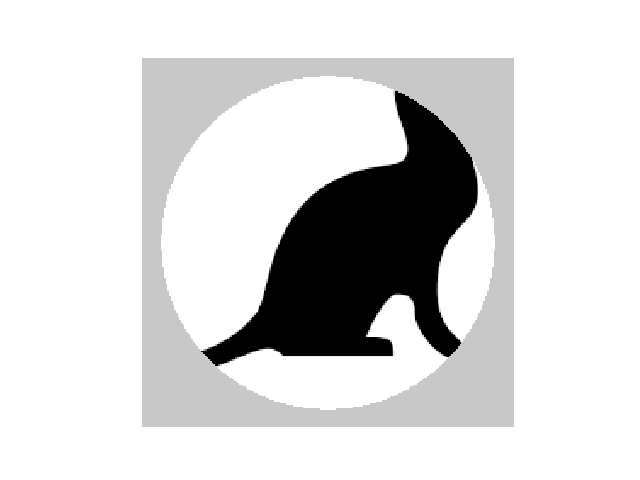}
			\caption{Applying a disk-shaped mask on the image of a cat.} \label{fig:cat_mask}
		\end{subfigure}
		\begin{subfigure}{.5\textwidth}
			\centering
			\includegraphics[scale=0.4]{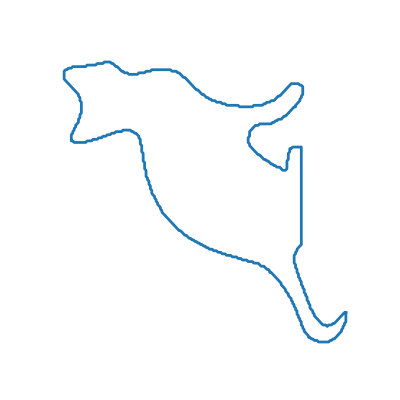}	
			\caption{Rotation of the contour image of the cat using the rotation matrix.} \label{fig:cat2}
		\end{subfigure}
		\mbe
		\caption{Application of mathematics in image processing using Scikit image processing library (Python).}
	\end{figure}

	\item  \textbf{Image masking.} Image masking is a method wherein a set of pixels are selected so that further image-manipulations can be performed on them. The Python code in Fig.~$ \ref{fig:cat_code_mask} $ is used to create a disk-shaped mask using the image of the cat shown in Fig.~$ \ref{fig:cat}$. The result of applying this mask is shown in Fig.~$ \ref{fig:cat_mask}$. 
	\begin{tcolorbox}
		The line of code which is used to create this disk shaped mask is 
		\begin{verbatim}
		outer_disk_mask = (X - l_x / 2)**2 +  
		(Y - l_y / 2)**2 > (0.45*l_x)**2.
		\end{verbatim}
		In fact,  all the pixels lying outside the disk of radius $ 0.45 l_x $ centered at $ \left(\frac{l_x}{2}, \frac{l_y}{2} \right) $ are getting selected. A student with knowledge of coordinate geometry can appreciate the fact that how mathematically this region is described by $ \displaystyle{\left(x-\frac{l_x}{2}\right)^2 + \left(y-\frac{l_y}{2}\right)^2 > (0.45*l_x)^2} $.  
	\end{tcolorbox}

\end{enumerate}

\begin{tcolorbox}
{
	Unfortunately, to the best of our knowledge at this moment, there is no image processing library available in the spirit of Turtle graphics (except perhaps the ``novice'' module of the Python Scikit Image available at \url {https://scikit-image.org/} to a certain extent), which is written for the pedagogical purpose of teaching image-processing-related programming to kids. However, working in the spirit of the above examples, such a library can be created, wherein one essentially hides the unnecessary program details from beginner students and still allow them to concentrate on the relevant mathematical concepts and their applications in image processing. Of course, we have just given a glimpse of  the possibilities with only a few examples. However, even with these few examples some of the advantages of this approach are obvious. 
	\begin{enumerate}
		\item It is a hands-on approach that allows a student to practically experiment with the code, 
		\item it allows students to see concrete real word applications of abstract mathematical concepts,
		\item it offers visual feedback and encourages students to express their creativity.
	\end{enumerate}
	Indeed, a student math art photo or math video-animation competition or exhibition can be organized in the classroom, enriching students' imagination and creativity, while they  also learn some interesting mathematical concepts during the process.  } 
\end{tcolorbox}

\subsection{{Physics and engineering.}}\label{Sect:Physics_Eng}
Mathematics is extensively used in Physics and Engineering. Indeed, it appears that the laws of nature are written in the language of mathematics. Unfortunately, often learning  mathematics is totally disconnected from learning physics. Renowned mathematician V.I. Arnold had once remarked (in an address on the teaching of mathematics in Paris in 1997) that ``In the middle of the twentieth century it was attempted to divide physics and mathematics. The consequences turned out to be catastrophic. Whole generations of mathematicians grew up without
knowing half of their science and, of course, in total ignorance of any other sciences. They first began teaching their ugly scholastic pseudo­mathematics to their students, then to schoolchildren (forgetting Hardy's warning that ugly mathematics has no permanent place under the Sun).'' We believe that Physics can be used to shape and inspire the learning of 
mathematics of a gifted child.

In the following, we give some examples and possible ideas for using concepts from physics and engineering in the math education of a gifted child.

\begin{enumerate}[label=Ex.~\arabic*.]
	\item \textbf{Electrical circuits and a positive definite matrix.} 
	A gifted student in high school, who has learned about positive definite matrices, can benefit from the example that we consider now. We say that a symmetric $ {\displaystyle n\times n} $ real matrix $ {\displaystyle A} $  is positive definite if the scalar $ {\displaystyle x^{\textsf {T}}Ax } $ is positive for every non-zero column vector $ {\displaystyle x} \in \mathbb{R}^n $.
	
	Suppose the matrix 
	\begin{equation} \label{Eq:defR}
	R = \mat{R_1 + R_2 + R_4}{R_2}{R_2}{R_2+R_3+R_5}
	\end{equation}
	is given such that $ R_1,R_2,R_3,R_4,R_5 >0 $. 
	
	We can use basic physics and the electrical circuit in Fig.~$ \ref{fig:circuit} $,  to show that the matrix $ R $ is positive definite. Using Kirchhoff's current law (KCL) and Kirchhoff's voltage law (KVL) we can write the following equation for the circuit. 
	\begin{align}
	\left[ \begin{smallmatrix}
	V_1
	\\ 
	V_2
	\end{smallmatrix}   \right] = \mat{R_1 + R_2 + R_4}{R_2}{R_2}{R_2+R_3+R_5} \left[ \begin{smallmatrix}
	i_1
	\\ 
	i_2
	\end{smallmatrix}   \right].
	\end{align}
	Let $ V = \left[ \begin{smallmatrix}
	V_1
	\\ 
	V_2
	\end{smallmatrix}\right] $ and $ I = \left[ \begin{smallmatrix}
	i_1
	\\ 
	i_2
	\end{smallmatrix}   \right]. $ 
	
	Since the electric power in a purely resistive circuit is always consumed by the load, we get that the total power delivered by the voltage sources is positive, i.e., $ V_1 i_1 + V_2 i_2 > 0 $ unless $ i_1=i_2=0 $. It is easy to see that,
	\begin{align}
	V_1 i_1 + V_2 i_2 > 0 \implies V^T I > 0 \implies (RI)^T I > 0 \implies I^T R^T I > 0 \implies I^T R I > 0, 
	\end{align}
	unless $ i_1=i_2=0 $. This is precisely the mathematical definition of a positive definite matrix mentioned earlier.
	We also note that in this example one can directly show without using the electrical circuit that the matrix $ R $ is positive definite.  An easy calculation shows that 
	\begin{align}
	I^T R I &= (R_1 + R_2 + R_4) i_1^2 + 2 R_2 i_1 i_2 + (R_2+R_3+R_5 )i_2^2  \nonumber\\
	&= (R_1 + R_4) i_1^2 + R_2 (i_1 + i_2)^2 + (R_3+R_5 )i_2^2 >0
	\end{align}
	unless $ i_1,i_2 =0 $.
	Therefore, the example illustrates the meaning of a matrix being positive definite in the context of the power dissipated in an electrical circuit. 
	\begin{figure}
		\centering
		\includegraphics[scale=0.8]{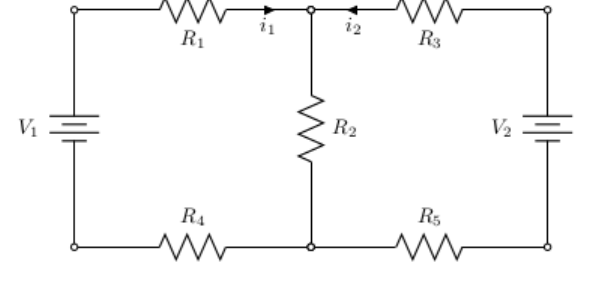}
		\caption{Application of mathematics in basic electrical circuit analysis. The above circuit can be used to show that the matrix $ R $, as given in Eq.~$\ref{Eq:defR} $, is positive definite.}\label{fig:circuit}
	\end{figure}
	\item \textbf{Electrical circuit voltage and current measurements.} Simple electrical circuits with a voltage source and a few resistors and possibly some LEDs are easy to analyze using  Kirchhoff's current and voltage laws. These circuits will result in solving a system of linear equations. Solutions to these equations can be checked using direct measurements by a multimeter and the differences in practical measurements and theoretical predictions can be analyzed. This is a good practical experiment that is accessible to gifted students in middle and high school. Moreover, the solution of linear equations leads to matrices and the rich subject of linear algebra.   
	\item \textbf{Bufoon's needle.} This is a well-known method of calculating the approximate value of $ \pi $ by a physical experiment of dropping a needle on a grid of parallel lines, whose spacing is greater than the length of a needle. The probability of the needle lying across a line is related to the value of $ \pi $. Suppose in an experiment out of $ n $ needles $ c $ of those needles crossed lines, then	it turns out that $ {\displaystyle \pi \approx {\frac {2l\cdot n}{t c}}}$, where  
	$ l $ is the length of the needle and $ t $ is the spacing between the parallel lines. Although, the proof of this fact involves calculus, for a bright high school student this experiment could be a good opportunity to learn about the concept of geometric probability. Moreover, approximating areas using many small rectangles can be introduced (basically, leading to the concept of Riemann sums).
	\item \textbf{Moving a robot on a given curve.} This is a slightly advanced challenge for a gifted student in high school who is interested in electronics and robotics. There are open-source single-board microcontroller kits, such as Arduino, available for building digital devices and interactive objects. It will be a challenging hands-on project to move a robot on a given curve, say a parabola or an ellipse. 
\end{enumerate}

\subsection{Computer games.}\label{Sect:Computer_games}
Computer games can be very effective learning tools. Games are hands-on and often involve motor control coupled with the processing of visual information. This helps in the understanding of concepts and also in long-term memory retention. 
%

In the following, we give some examples of games that can be helpful in teaching a gifted student. 

\begin{enumerate}[label=Ex.~\arabic*.]
	\item \textbf{Negative numbers on a ladder and a monkey.} A simple game in which a monkey can jump up or down on a ladder (see Fig.~$ \ref{fig:ladder} $). The up direction is positive and the down direction is negative. Suppose a monkey is at position $ x $ and there is a banana at position $ y $ on the ladder. Assume that the monkey can move one step on the ladder in one jump. The player is required to enter the correct number of jumps with the correct sign to move the monkey to the banana. This game can be used to introduce negative numbers. Moreover, some variations of this game, such as changing the jump from one step to two, three, or more steps can help in reinforcing the multiplication facts.
	
	\begin{figure}
		\centering 
		\includegraphics[]{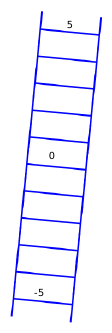}
		\caption{A simple game in which a monkey can jump up or down on a ladder. This game can be used to teach `negative numbers' to students.} \label{fig:ladder}
	\end{figure} 
	\item \textbf{Shooting games to learn coordinate geometry.} A target is required to hit and for that the coordinate of the target is needed. This game can have several variations, such as when the target is one vertex of a right-angle triangle and the player is required to use the Pythagorean theorem to find the coordinate of the target. The game setup can include objects such as walls, roofs, catenary arches, flowers ect., providing a blending of real-life objects with nice mathematical descriptions. \label{enum:shootinggame}
	
	\item  \textbf{Calculus games.} Perhaps it would be somewhat more challenging to design games for learning calculus. Still, there are many possibilities. For example, a game could be set up in a desert where the goal of the game could be to design a tent using limited resources. This tent design problem could be such that an optimization problem is required to be solved in order to win. Another scenario might require some approximate calculation based on the concept of the Riemann sum.   An engaging interactive game, with a series of tasks involving calculus concepts, might be very helpful for students learning calculus.
\end{enumerate}

\section{Dealing with failures.}\label{Sect:Failure}
So far we have seen some strategies to engage gifted kids and get them interested in mathematics. Now, we discuss an important topic, which is dealing with failures while learning mathematics. Unfortunately, this topic is often not given its due importance in a typical math classroom. While teaching a class of young students, the teacher asked a simple addition question. Many students raised their hands to answer the question. The teacher picked Danny and asked him to answer. Danny gave a wrong answer and as soon as the teacher said no it was wrong, Danny started crying, and with tearful eyes, he swiftly ran away from the class. Danny is an extreme example, but it is still true that many students run away from mathematics because they can not tolerate failure. Negative emotions associated with failure in mathematics, such as failure to answer a question in class or low grades in the subject begin to affect the self-worth of an individual learner. As a coping mechanism, they start to avoid mathematics as much as possible. This exacerbates the situation, and such students further fall behind in their mathematical learning goals. Later on, many of such students become adults who are intelligent and successful in their professions, but they still resent mathematics and proudly proclaim that -``I'm not a math-person''. One remedy to this lies in the hands of teachers teaching young students. The teachers should emphasize that failure in mathematics is inevitable. Even professional research mathematicians struggle in solving their problems. There are many open problems, some for several decades. Therefore it is important to enjoy the process of learning without any worry of failure. This message should be conveyed to every young student of mathematics.

\bibliographystyle{apalike}

\end{document}